\documentclass[a4paper]{article}
\usepackage[a4paper,left=3cm,right=2cm,top=2.5cm,bottom=2.5cm]{geometry}
\usepackage{graphicx} 
\usepackage{amsmath,amsthm}
\usepackage{epstopdf} 
\usepackage{mathtools}
\usepackage{amsfonts}
\usepackage{amssymb}
\usepackage{hyperref}
\usepackage{xcolor}
\usepackage{matlab-prettifier}
\usepackage[normalem]{ulem}
\usepackage{comment}
\usepackage[algoruled]{algorithm2e}
\usepackage{mathabx} 
\usepackage{subcaption}
\usepackage{booktabs}

\usepackage[shortlabels]{enumitem}

\newtheorem{theorem}{Theorem}[section]
\newtheorem{proposition}[theorem]{Proposition}

\newtheorem{lemma}[theorem]{Lemma}
\newtheorem{corollary}[theorem]{Corollary}
\theoremstyle{definition}

\newtheorem{remark}[theorem]{{\sc Remark}}
\newtheorem{example}[theorem]{Example}

\newtheorem{assumption}[theorem]{Assumption}

\def\R{\mathbb{R}}
\def\C{\mathbb{C}}

\def\F{\mathbb{F}}

\DeclareMathOperator{\rank}{rank}
\DeclareMathOperator{\im}{Im}

\DeclareMathOperator{\tr}{trace}

\DeclareMathOperator{\argmin}{argmin}
\DeclareMathOperator{\vvec}{vec}
\DeclareMathOperator{\spann}{span}

\usepackage{tikz, pgfplots, mathtools}
\pgfplotsset{compat=1.18}

\DeclarePairedDelimiter{\norm}{\lVert}{\rVert}
\DeclarePairedDelimiter{\abs}{\lvert}{\rvert}

\newcommand{\subjectto}{\ensuremath{\ \text{s.t.}\ }}

\definecolor{brilliantrose}{rgb}{1.0, 0.33, 0.64}
\definecolor{myviolet}{rgb}{0.21, 0.0, 0.85}
\definecolor{amethyst}{rgb}{0.6, 0.4, 0.8}
\definecolor{carrotorange}{rgb}{0.93, 0.57, 0.13}

\title{Riemann-Oracle: \\
A general-purpose Riemannian optimizer to solve nearness problems in matrix theory}

\author{Miryam Gnazzo\thanks{
University of Pisa, Department of Mathematics, I-56127 Pisa, Italy. Affiliated to the Italian INdAM-GNS (Gruppo Nazionale di Calcolo Scientifico). Supported 
by the PRIN 2022 Project \emph{Low-rank Structures and Numerical Methods in Matrix and Tensor Computations and their Application}. Email: miryam.gnazzo@dm.unipi.it}, Vanni Noferini\thanks{Corresponding author. Aalto University, Department of Mathematics and Systems Analysis, P.O. Box 11100, FI-00076, Aalto, Finland. Supported from a Research Council of Finland grant (decision number 370932), which has been awarded and will begin on 1 September 2025. Email: vanni.noferini@aalto.fi}, Lauri Nyman\thanks{Aalto University, Department of Mathematics and Systems Analysis, P.O. Box 11100, FI-00076, Aalto, Finland. Email: lauri.s.nyman@aalto.fi}, and Federico Poloni\thanks{University of Pisa, Department of Computer Science, I-56127 Pisa, Italy. Affiliated to the Italian INdAM-GNCS (Gruppo Nazionale di Calcolo Scientifico). Supported by INdAM-GNCS 2024 research project CUP\_E53C23001670001, and by the European Union-NextGenerationEU National Recovery and Resilience Plan under the PRIN 2022 grant \emph{Low-rank Structures and Numerical Methods in Matrix and Tensor Computations and their Applications} 20227PCCKZ – CUP\_I53D23002280006 and under the National Centre for HPC, Big Data and Quantum Computing / Spoke 6, \emph{Multiscale Modeling and Engineering Applications}. Email: federico.poloni@unipi.it}
}
\date{\today}

\begin{document}
\maketitle

\begin{abstract}
    We propose a {versatile} approach to address a large family of matrix nearness problems, possibly with additional linear constraints. Our method is based on splitting a matrix nearness problem into two nested optimization problems, of which the inner one can be solved either exactly or cheaply, while the outer one can be recast as an unconstrained optimization task over a smooth real Riemannian manifold. We observe that this paradigm applies to many matrix nearness problems of practical interest appearing in the literature, thus revealing that they are equivalent in this sense to a Riemannian optimization problem. We also show that the objective function to be minimized on the Riemannian manifold can be discontinuous, thus requiring regularization techniques, and we give conditions for this to happen. Finally, we demonstrate the practical applicability of our method by implementing it for a number of matrix nearness problems that are relevant for applications and are currently considered very demanding in practice. Extensive numerical experiments demonstrate that our method often greatly outperforms its predecessors, including algorithms specifically designed for those particular problems.
\end{abstract}


\bigskip


\bigskip

\section{Introduction}
Nearness problems are a fundamental class of problems in numerical linear algebra and matrix theory. A classical survey of some important applications appeared in \cite{nicksurvey}; recent high-impact work on specific instances of nearness problems, or on the subject in general, includes for example \cite{Bini2010, nearsingpen, dhillontropp, nickspd, nickcorr, qisun, VL85}. Below, we provide a general description of matrix nearness problems, following in part the approach proposed in Nick Higham's PhD thesis \cite{nickthesis}.

 Given an input matrix $A$ and a property $\mathfrak{P}$ that does not hold for $A$, the corresponding matrix nearness problem consists of finding a matrix $B$ nearest to $A$ and such that $\mathfrak{P}$ holds for $B$, as well as the distance between $A$ and $B$. More formally, denoting by $\mathcal{Q}$ the set of matrices that have the property $\mathfrak{P}$, we seek the minimum and an argument minimum for the constrained optimization problem
\begin{equation}\label{eq:problem}
    \min_{{X\in \mathcal{Q}}} \| A - X \|. 
\end{equation} 
In \eqref{eq:problem}, we assume to have fixed a matrix norm. The most frequent choice in practice is the Frobenius norm, ${\| X \|_F = \sqrt{\tr X^* X}}$; while other choices can sometimes be of interest \cite{bora,dhillontropp}, in this paper we will focus on nearness problems with respect to the Frobenius distance. Matrix nearness problems can also be generalized to \emph{matrix pencils}, or more generally \emph{matrix polynomials}. In this case $A(z)=\sum_{i=0}^d A_i z^i$ is a polynomial with matrix coefficients, and one looks for a second matrix polynomial $B(z)=\sum_{i=0}^d B_i z^i$ that (a) minimizes the squared distance $\sum_{i=0}^d \| A_i - B_i \|^2$, and (b) has a given property $\mathfrak{P}$ that $A(z)$ lacks. This framework includes nearness problems for scalar polynomials, by interpreting their coefficients as $1\times 1$ matrices.

Some matrix nearness problems are fully understood in the sense that either they admit a closed-form exact solution, or a global minimizer can be computed numerically via some well-known (and fast) algorithm. Let us give one example for each situation. It is an elementary exercise to show that the matrix nearest to $A$ and such that its $(1,1)$ element is zero is $B=A-e_1 e_1^T A e_1 e_1^T$, with distance $|e_1^T A e_1|$; whereas it is a corollary of the Eckart–Young–Mirsky theorem \cite{EY36} that the singular matrix nearest to $A \in GL(n,\C)$ can be obtained by computing an SVD $A=\sum_{{i=1}}^n \sigma_i u_i v_i^*$ and then setting $B=\sum_{{i=1}}^{n-1} \sigma_i u_i v_i^*$, with distance $\sigma_n$. {On the other hand, other matrix nearness problems are well known to be hard \cite{DNN24, NP21}, and existing algorithms cannot promise anything more than a local minimum.} Often, current algorithms are also quite slow in practice.

 In this paper, we propose an approach that
 \begin{enumerate}
     \item is {versatile} in the sense that it can deal with a very large class of matrix nearness problems, including many of the instances that we could find described in the literature;
     \item is {competitive} in practice, to the point of often outperforming, in a vast majority of our experiments, currently existing algorithms for ``difficult" matrix nearness problems (and in spite of the fact that the best competitor algorithms are often specialized to one particular nearness problem, while our approach is very general); this is demonstrated by numerical experiments in Section \ref{sec:numericalexperiments}.
     \item ultimately relies on recasting a matrix nearness problem as an optimization task over a real Riemannian manifold.
 \end{enumerate}
 
 To illustrate the last item, recall that, in theoretical computer science, an \emph{oracle} is an abstract Turing machine that has access to a black-box with the ability to solve, in a single operation, any instance of a certain class of problems \cite{papa,sipser2006}. Our method capitalizes on the insight that many matrix nearness problems become more tractable when provided by an oracle with specific supplementary information about the minimizer. For instance, if we are working over square matrices, this information could be a certain optimal eigenvalue or eigenvector, say, $\theta$. Suppose further that, by restricting to matrices that share this information $\theta$, the problem~\eqref{eq:problem} can be solved exactly. In other words, let us write 
 the feasible set as $\mathcal{Q}=\bigcup_\theta \mathcal{Q}_\theta$, where $\mathcal{Q}_\theta$ is the set of matrices having both the property $\mathfrak{P}$ and the attribute $\theta$. Assume that we can compute the solution $f(\theta)$ of the stiffened problem
 \begin{equation}\label{eq:stiffendproblem}
    f(\theta) = \min_{{X \in \mathcal{Q}_\theta}} \| A - X \|.
 \end{equation}
At this point, we can solve the original problem \eqref{eq:problem} by optimizing \eqref{eq:stiffendproblem} over $\theta$, i.e., we can recast \eqref{eq:problem} as the equivalent optimization task
 \[ \min_\theta f(\theta).\]
Sometimes, in the literature, the functional $f$ to be minimized in \eqref{eq:stiffendproblem} is called the \textit{variable projection functional} \cite{variableprojection,oleary}. This minimization often has to be performed over an appropriate matrix or vector manifold; going back to the concrete example where $\theta$ is an eigenvector, and there are no additional restrictions on it, then one can optimize over the set of all possible normalized eigenvectors, i.e., the unit sphere. Or possibly, the information $\theta$ could consist of an arbitrary set of $d$ linearly independent eigenvectors; in this case, one would then optimize over the Grassmannian of $d$-dimensional subspaces \cite{EdelmanGrassmann}.

For some specific matrix nearness problems, optimization over manifolds has been already used as a tool, although usually in a more direct way that in our manuscript. We mention a few illustrative recent examples: Vandereycken \cite{bartvander} solved the low-rank matrix completion problem -- recovering a low-rank matrix from partially observed data -- by means of optimization over the manifold of matrices of fixed rank. More recently,  Goyers, Cartis, and Eftekhari \cite{GCE23} also addressed low-rank completion via Riemannian optimization, using a reformulation over the Grassmann manifold. Borsdorf \cite{borsdorf} applied an augmented Lagrangian algorithm to tackle a nearness problem motivated by chemistry and directly expressible as a minimum over the Stiefel manifold. On the other hand, there have been successful attempts to use structured total least norm algorithms \cite{RPG} and regularization for specific nearness problems. We provide but one recent example: Kaltofen, Yang and Zhi \cite{KYZ} {used this idea} in the context of approximate GCD; however, unlike what we do in the present paper, neither \cite{KYZ} nor related papers on GCD computations seem to have combined these techniques with Riemannian optimization.

The ideas proposed in this paper are applicable to a wide class of nearness problems. In essence, they amount to a two-level optimization framework, in which the inner subproblem can be solved cheaply and/or explicitly, while Riemannian optimization becomes a crucial tool for the outer subproblem. 

The roots of this ``stiffen-and-optimize" approach can be found in the literature, at least for some specific nearness problems. An early archetype, due to Byers and Van Loan, dates back to the 1980s~\cite{Bye88,VL85}: Given an $n \times n$ matrix $A$ with all its eigenvalues in the open left half-plane, one would like to find its distance to Hurwitz instability, i.e., the smallest $d = \norm{A-B}_F$ such that $B$ has an eigenvalue in the closed right-half plane. By continuity, this eigenvalue $\lambda = i\omega$ must lie on the imaginary axis. If an oracle gives us the value of $\omega \in\mathbb{R}$, then this problem reduces to finding the distance from singularity of $A-i\omega I$. As previously discussed, the solution is $d = \sigma_{n}(A- i\omega I)$. Hence, the original problem is equivalent to minimizing the function $f(\omega) = \sigma_{n}(A-i\omega I)$: This is a univariate optimization problem over $\omega \in\mathbb{R}$, and is simple to solve. The resulting algorithm is cheap and efficient, but it has the drawback that it cannot deal with additional structural constraints on the perturbation $B-A$. {Another early example of this approach for matrix nearness problems is \cite{KOM09}, where Keshavan, Oh and Montanari solved the low-rank matrix completion problem by optimizing exactly over the singular values of the solution, and then numerically over the singular vectors.}

In a series of papers~\cite{MarU13missing,MarU14software,UseM14manifold,UseM14,UseM17gcd}, Markovsky and Usevich solved structured low-rank approximation problems with an analogous two-level optimization approach: They imposed that the block matrix $\begin{bsmallmatrix}\Theta & -I\end{bsmallmatrix}$ is in the left kernel of a matrix $A$ with a prescribed structure, derive an expression for the objective function $f(\Theta)$ as the solution of a least-squares problem (similarly to what we shall do in Section~\ref{sec:matrix}), and then minimize it (using different techniques from the ones that we propose in this paper). Their focus is mainly on matrices with Hankel blocks, but the methods extend to more general settings.

More recently, Noferini and Poloni \cite{NP21} proposed a Riemannian algorithm for the problem of finding the \emph{nearest $\Omega$-stable matrix}, which is related to the computation of the distance to Hurwitz instability, but more challenging:
Given a square matrix $A$, one looks for the nearest matrix $B$ whose eigenvalues all lie in a prescribed closed set $\Omega$. In this case, the problem becomes simple if an oracle prescribes a unitary matrix $Q$ such that $Q^*BQ$ is in Schur form, and the method in \cite{NP21} optimizes over $Q$; however, this algorithm also seems difficult to extend to the case where additional constraints are added. A similar idea was used in \cite{DNN24} by Dopico, Noferini and Nyman for the problem of finding the nearest singular pencil \cite{nearsingpen}, again achieving remarkable success with respect to alternative approaches. 

Although not explicitly derived in this spirit, the method studied by Das and Bora in \cite{bora} to find the nearest singular polynomial matrix is also a special case of the more general paradigm proposed in this paper. Indeed, we will see in Section \ref{sec:Nearest_sing_matrix_poly} that our approach, when specialized to the nearest singular polynomial problem, leads to the same core formula underlying the algorithm in \cite{bora}; nonetheless, our practical implementation is substantially different than \cite{bora}, because we use Riemannian optimization and regularization of the objective function. Our numerical experiments, reported in Subsection \ref{sec:betterthanbora}, show that these differences lead to a significantly better performance with respect to \cite{bora}. This test case, together with various other experiments presented in Section \ref{sec:numericalexperiments}, provides convincing evidence that it is tremendously effective to tackle matrix nearness problem by means of the combined power of ``oracles", Riemannian optimization, and regularization techniques. We call our resulting general-purpose algorithm \emph{Riemann-Oracle}.

In the rest of the article, we describe both the theoretical and implementative aspects of Riemann-Oracle by going through several classes of nearness problems that can be tackled with our algorithm. Among the many nearness problems potentially solvable with Riemann-Oracle, we have singled out a handful of significant (in our opinion) examples. Our selection was based on two metaphorical metrics: The examples should be of current pragmatic interest (for example in numerical linear algebra, control theory, or computer algebra), and illustrative of the properties of our method. 

We start in Section \ref{sec:matrix} by studying the fundamental problem of computing the nearest singular matrix to a given one, while only allowing perturbations that belong to a prescribed vector space $\mathcal{S}$. We expose the theory behind Riemann-Oracle in detail for this case, and our analysis leads us to the task of minimizing a possibly discontinuous function over a smooth real manifold. This objective function can be optimized via classical regularization tools such as penalty methods or augmented Lagrangian algorithms. Section \ref{sec:riemann} discusses how the inner optimization problem in those general algorithms can be tackled over a real Riemannian manifold, and Section \ref{subsec:convergence-results} discusses convergence of the outer iteration. For specific vector spaces $\mathcal{S}$, the computation described in Section \ref{sec:matrix} can be simplified; we illustrate this in Section \ref{sec:sparse} for the case of matrices with a prescribed sparsity pattern. We consider as a pivotal test case the structured distance to singularity in Section \ref{sec:matrix}, as many other problems can be reduced to it. We illustrate this fact in Section \ref{sec:Nearest_sing_matrix_poly}, Section \ref{sec:approxgcd}, and Section \ref{sec:instability}. Therein we consider, respectively, the problems of finding the nearest singular polynomial matrix (given an upper bound on its degree), the approximate GCD (given a lower bound on its degree) of two scalar polynomials, and the nearest matrix with at least one eigenvalue in a prescribed closed set (distance to instability). In Section \ref{sec:nearest_matrix_pescribed_nullity}, we go beyond the structured distance to singularity and seek the nearest matrix having prescribed rank, while still restricting the set of possible perturbations to an arbitrary, but fixed, linear subspace; this requires changing the manifold over which optimization is performed, and showcases the flexibility of Riemann-Oracle. Section \ref{sec:numericalexperiments} describes extensive numerical tests in which we compare the performance of Riemann-Oracle to pre-existing algorithms in all the examples mentioned above.
A striking common theme that can be observed from the results of the experiments in Section \ref{sec:numericalexperiments} is that, in spite of its broader applicability spectrum, Riemann-Oracle is almost\footnote{The one experiment in which Riemann-Oracle performed slightly worse than (some) pre-existing algorithms was for a very challenging approximate GCD problem, see Subsection \ref{sec:herewelose}.} always competitive with, and in several cases it significantly outperforms, algorithms specialized on just one specific class of nearness problem in terms of computational speed, or quality of the answer, or both. Finally, we draw some conclusions in Section \ref{sec:conclusions}. A technical proof is deferred to Appendix \ref{sec:proofgradienthessian}.

\section{Nearest structured singular matrix} \label{sec:matrix}

We start our exposition by analyzing the problem of computing the singular matrix $B=A+\Delta$ nearest to a given matrix $A$, where the perturbation $\Delta$ must be optimally picked from a prescribed vector space $\mathcal{S}$. This is a very deliberate choice; indeed, as we shall see later on, many other matrix nearness problems can be reformulated as a nearest structured singular matrix problem.

Let the field $\mathbb{F}$ be either $\mathbb{R}$ or $\mathbb{C}$. We consider a matrix $A\in\mathbb{F}^{m\times n}$, with $m\geq n$, and we seek its possible perturbations in the linear span of $p$ given (linearly independent) matrices $P^{(1)}, \dots, P^{(p)} \in \mathbb{F}^{m\times n}$, i.e.,
\[
\mathcal{S} = \left\{ \Delta \in \mathbb{F}^{m\times n} \colon \Delta = \sum_{i=1}^p P^{(i)} \delta_i, \quad \delta = \begin{bmatrix}
    \delta_1\\
    \vdots\\
    \delta_p
\end{bmatrix} \in \mathbb{F}^p \right\}.
\]
Alternatively, we can write
\[
\vvec \Delta = \mathcal{P}\delta, \quad \mathcal{P} = \begin{bmatrix}
    \vvec P^{(1)} & \vvec P^{(2)} & \dots & \vvec P^{(p)}
\end{bmatrix} \in \mathbb{F}^{mn \times p}.
\]
\begin{remark}
    In even greater generality, we could consider a basis $\{ P^{(i)}\}$ over $\R$ even when $\F=\C$. Since, eventually, we are going to construct a real Riemannian manifold, it is no loss of generality for our next developments to assume that $\mathcal{S}$ is always a real vector space. However, assuming that $\delta \in \F^p$ (as opposed to $\delta \in \R^p$ even when $\F=\C$) does simplify some formulae, and since this is the most common case in practice, for the sake of exposition we make this choice.
\end{remark}

\begin{example}\label{ex:first} The set $\mathcal{S}$ of symmetric Toeplitz $3\times 3$ matrices is a $3$-dimensional vector space over $\mathbb{F}$. Indeed, it can be described as the set of linear combinations of
\[
P^{(1)} = \frac{1}{\sqrt{3}}\begin{bmatrix}
    1 \\ & 1\\ && 1
\end{bmatrix}, \quad
P^{(2)} = \frac{1}{2}\begin{bmatrix}
    & 1\\
    1 && 1\\
    & 1
\end{bmatrix}, \quad
P^{(3)} = \frac{1}{\sqrt{2}}
\begin{bmatrix}
    &&1\\
    \\
    1
\end{bmatrix},
\]
since every symmetric Toeplitz matrix $\Delta$ can be written as $P^{(1)}\delta_1 + P^{(2)}\delta_2 + P^{(3)}\delta_3$, for a suitable choice of the vector $\delta\in\mathbb{F}^3$.
\end{example}
In Example \ref{ex:first}, the $P^{(i)}$ are chosen so that they are orthogonal, i.e., $\tr [P^{(i)}]^* P^{(j)}=0$ when $i \neq j$, and their Frobenius norm is $1$. More generally, we make the following assumption throughout the paper.
\begin{assumption}\label{assumption22}
We assume in the following that $\mathcal{P}^*\mathcal{P} = I_p$, i.e., the columns of $\mathcal{P}$ are orthonormal. Hence $\norm{\Delta}_F = \norm{\delta}$.
\end{assumption}
Assumption \ref{assumption22} holds naturally in many cases, for instance when the nonzero elements of the matrices $P^{(i)}$ are in disjoint locations. Moreover, it does not cause any loss of generality, because if $\mathcal{P}$ is a generic matrix then we can obtain an equivalent perturbation basis that satisfies the assumptions, for example by computing a (rank-revealing) thin QR factorization and replacing $\mathcal{P}=QR$ and $\delta$ with $Q$ and $R \delta$.

Note that $A+\Delta$ is singular if and only if there exists a vector $v$ satisfying both $\norm{v}=1$ and $(A+\Delta)v = 0$. We next show in Theorem \ref{thm:fv} that, if an oracle were to provide the correct null vector $v$ for the optimal $\Delta$, then we could solve exactly the problem of finding $\Delta$.  In other words, we claim that it is easy to evaluate
\begin{equation} \label{fv}
    f(v) = \min_{{\Delta\in\mathcal{S}}} \, \norm{\Delta}_F^2 \subjectto (A+\Delta)v = 0.
\end{equation}

\begin{theorem} \label{thm:fv}
Let a vector $v\in\mathbb{F}^n$, $v\neq 0$, be given. Define
\begin{equation} \label{Mr}
    M = M(v) = \begin{bmatrix}
        P^{(1)}v & P^{(2)}v & \dots & P^{(p)}v
    \end{bmatrix} \in \mathbb{F}^{m \times p}, \quad r = r(v) = -Av.
\end{equation}
Then:
\begin{enumerate}
    \item The minimization problem in~\eqref{fv} is well-posed (i.e., the minimum is sought over a non-empty set) if and only if $r \in \operatorname{Im} M$. When this happens, the minimum in~\eqref{fv} is attained for 
\[
    \Delta_* = \sum_{i=1}^p P^{(i)} (\delta_*)_i, \quad  \delta_* = M^\dagger r,
\]
where $M^\dagger$ denotes the Moore--Penrose pseudoinverse.

If, in addition, $M$ has full row rank, then we have explicit expressions
\[
\delta_* = M^*(MM^*)^{-1}r, \quad f(v) = r^*(MM^*)^{-1}r.
\]
\item Defining the least-squares function $g(v) := \| M(v)^\dagger r(v) \|_F^2$ for all $v \in \F^n$, then
\[  f(v) = \begin{cases}
    g(v) \ &\mathrm{if} \ r(v) \in \im M(v);\\
    +\infty \ &\mathrm{otherwise}.
\end{cases}     \]
\end{enumerate}
\end{theorem}
\begin{proof}
\begin{enumerate}
    \item  Since $\Delta v = M\delta$ for each $\delta$, we can rewrite the condition $(A+\Delta)v = 0$ as $M\delta = r$. We know already that $\norm{\Delta}_F^2=\norm{\delta}^2$, and hence
 \[
 f(v) = \min_{{\delta \in \mathbb{F}^p}}\, \norm{\delta}^2 \subjectto M\delta = r.
 \]
By the standard theory of least-squares problems (see, e.g.,~\cite[Theorem~1.2.10]{Bjorck}), the solution can be obtained by left-multiplying the right hand side with the Moore--Penrose pseudoinverse
\[
\delta_* = M^\dagger r = \arg\min_{{\delta \in \mathbb{F}^p}} \,\norm{\delta}^2 \subjectto M\delta = r.
\]
When $M$ has full row rank, the Moore--Penrose pseudoinverse has the explicit formula $M^\dagger = M^*(MM^*)^{-1}r$, and we see that
\[
f(v) = \norm{\delta_*}^2 = (\delta_*)^*\delta_* = r^*(MM^*)^{-1}MM^*(MM^*)^{-1}r = r^*(MM^*)^{-1}r. 
\]
\item The previous item immediately yields $f(v)=g(v)$ when $v$ is such that $r \in \im M$. On the other hand, if $v$ lies outside this feasible region, then for such $v$ the set of perturbations $\Delta$ over which the minimum is sought in \eqref{fv} is empty, and by definition ${\min{}} \emptyset = +\infty$.
\end{enumerate}
\end{proof}
An important special case happens when $A \in \mathcal{S}$, i.e., we are looking for a perturbation with a structure that is shared by $A$; for instance, Toeplitz perturbations of a Toeplitz matrix, or sparse perturbations of a sparse matrix (with the same zero pattern). In this case, the solution takes the form of an orthogonal projection. Note that, in the statement of Corollary \ref{cor:fv}, $M$ depends on $v$ but $\alpha$ does not. 
\begin{corollary} \label{cor:fv}
If $A = \sum_{i=1}^p P^{(i)}\alpha_i$ for a certain vector $\alpha \in \mathbb{F}^p$, then $Av = M\alpha$, and
    \[
    \delta_* = -M^\dagger M\alpha, \quad
    \alpha+\delta_* = (I - M^\dagger M)\alpha,
    \]
    i.e., $-\delta_*$ is obtained by projecting $\alpha$ orthogonally onto $\operatorname{Im} M^*$, and $\alpha+\delta_*$ by projecting it onto $(\operatorname{Im} M^*)^\perp$.  
\end{corollary}

Observe that $f(v)$ is $0$-homogeneous, i.e., $f(v\beta ) = f(v)$ for each $\beta \in\mathbb{F}$, $\beta \neq 0$. This suggests computing the minimum over the set of norm-1 vectors, i.e., the unit sphere in $\mathbb{F}^n$, henceforth denoted by $\mathcal{M}$.

Thus, we have
\[
    \min_{{\Delta\in\mathcal{S}}} \norm{\Delta}_F^{{2}} \subjectto \text{$A+\Delta$ is singular} = \min_{{v \in \mathcal{M}}} f(v).
\]
The equality is clear, since $A+\Delta$ is singular if and only if there exists a unit-norm vector $v$ such that $(A+\Delta) v = 0$.

\begin{remark}
    In some applications, one might desire to study the variant of the nearest structured singular matrix problem where $A \not \in \mathcal{S},\Delta \not \in \mathcal{S}$, but $A+\Delta \in \mathcal{S}$. While we did not directly analyze this case, it is easy to see that it is equivalent to the case $A \in \mathcal{S},\Delta \in \mathcal{S}$, studied above, by a linear change of variables. Indeed, consider the (unique) decomposition $A=A_\mathcal{S} + A_\perp$ where $A_\mathcal{S} \in \mathcal{S}, A_\perp \in \mathcal{S}^\perp$. Then, for all $\Delta$ such that $A+\Delta \in \mathcal{S}$ we have, defining $\Delta_\mathcal{S}:=A + \Delta - A_\mathcal{S} \in \mathcal{S}$,
    \[ \| \Delta \|^2_F = \| \Delta_\mathcal{S} \|_F^2 + \|A_\perp \|_F^2  \]
    and hence the problem of finding the singular matrix $A+\Delta \in \mathcal{S}$ nearest to $A$ is equivalent to finding the singular matrix $A_\mathcal{S}+\Delta_\mathcal{S}$ nearest to $A_\mathcal{S}$.
\end{remark}

The full-rank case of Theorem~\ref{thm:fv} and the expression $f(v) = r(v)^* (M(v)M(v)^*)^{-1}r(v)$ already appear in an equivalent form in the work of Markovsky and Usevich; see e.g.~\cite[Section~3.1]{UseM14}. However, in the next section we will see that the full-rank assumption can be too restrictive for certain problems.

\subsection{The need for regularization}
\label{sec:need_regularization}
Unfortunately, the assumption that $M$ has full row rank is restrictive in many practical cases. We give an illustrative example below.

\begin{figure}
    \centering
    \includegraphics[width=0.6\textwidth]{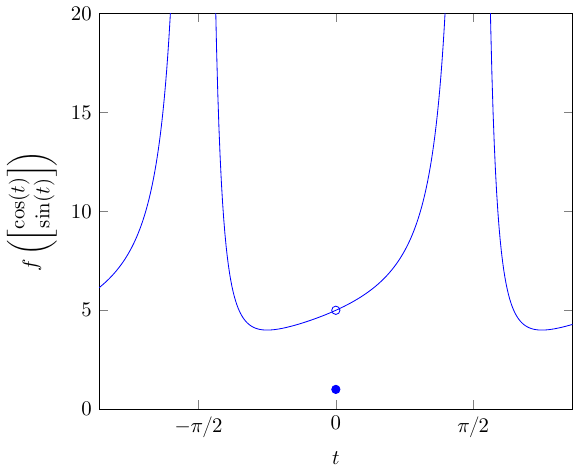}\\[1em]
    \includegraphics[width=0.75\textwidth]{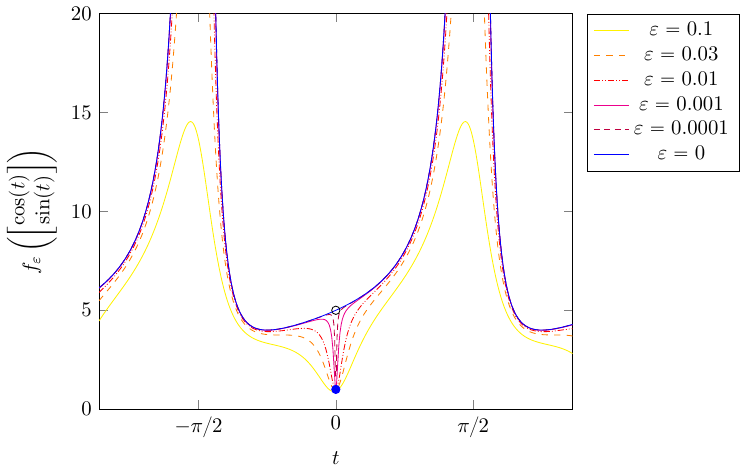}
    \caption{Top: Plot of $f(v)$ as in \eqref{fv} for Example \ref{ex:needy}. Bottom: Plot of $f_\varepsilon(v)$ as in \eqref{fepsilonv}, for selected values of $\varepsilon>0$,  for Example \eqref{ex:needy}. The label $\varepsilon=0$ indicates $f(v)$.}
    \label{fig1}
\end{figure}

\begin{example}\label{ex:needy}
    Over $\F=\R$, consider the matrix $A=\begin{bsmallmatrix}
        1&1\\
        0&2
    \end{bsmallmatrix}$. We look for real perturbations of the diagonal elements only, that is, $P^{(1)}=\begin{bsmallmatrix}
        1&0\\
        0&0
    \end{bsmallmatrix}$, $P^{(2)}=\begin{bsmallmatrix}
        0&0\\
        0&1
    \end{bsmallmatrix}$, and $A \not\in \mathcal{S}$. For a given unit vector $v=\begin{bsmallmatrix}
        \cos(t)\\
        \sin(t)
    \end{bsmallmatrix} \in \mathcal{M}$, parametrized by $t \in \R$, we obtain
    \[  M = \begin{bmatrix}
        \cos(t) & 0\\
        0 &\sin(t)
    \end{bmatrix}, \qquad r=\begin{bmatrix}
    -\sin(t)-\cos(t)\\
        -2\sin(t)
    \end{bmatrix} .   \]
    Hence  (see also Figure \ref{fig1}),
    \[  f(v(t)) = \begin{cases}
      4 + (1 + \tan(t))^2  &\mathrm{if} \ 2t \neq k \pi \ (k \in \mathbb{Z});\\
      1  &\mathrm{if} \ t = k \pi \ (k \in \mathbb{Z});\\
      +\infty  &\mathrm{if}  \  t = \frac{k \pi}2 (k \ \mathrm{odd}).\\
    \end{cases}   \]
    It is clear that in this example $f(v)$ is discontinuous, even if we restrict it to the feasible region on which $f(v)=g(v)$, i.e., we exclude the points $\pm e_2$. In particular $f(\pm e_1)=1$, but $\lim_{v \rightarrow \pm e_1} f(v)=5$; note that $1$ is also the global minimum of $f(v)$, thus demonstrating that a discontinuity can possibly happen at a minimizer.
\end{example}    
    
In Example \ref{ex:needy}, $A \not \in \mathcal{S}$. If $A \in \mathcal{S}$ then $f(v)$ is bounded on the unit sphere, since by Theorem \ref{thm:fv} and Corollary \ref{cor:fv},
\[ \min_{{v\in\mathcal{M}}} f(v) = \min_{{v\in\mathcal{M}}} g(v) \leq \min_{{v\in\mathcal{M}}} \|M(v)^\dagger M(v)\| \| \alpha \|=\| \alpha \|.  \] 
However, even when $A \in \mathcal{S}$, the function $f(v)$ may have jump discontinuities, even at a global minimizer. This can be seen for instance by slightly modifying Example \ref{ex:needy} to $A=\begin{bsmallmatrix}
    1&0\\
    0&2
\end{bsmallmatrix}$. If one employs optimization algorithms designed for smooth functions, the task of minimizing numerically a discontinuous function is (almost) hopeless. A standard remedy to this issue, and especially popular in the context of least-squares problems, is to employ Tikhonov regularization. Figure \ref{fig1} also illustrates the effect of regularization to Example \ref{ex:needy}; in the next subsection, we analyze the general situation.

\begin{remark} \label{rem:smallcounterex}
    The approach of~\cite{MarU14software, UseM14} is similar to ours, but it assumes that $M$ has full rank, and hence it is vulnerable to the issues just described. This can be revealed by slightly modifying the example above considering $A = \begin{bsmallmatrix}2 & 1\\0 & 1\end{bsmallmatrix}$ (the reasons for these minor modifications are that \cite{MarU14software, UseM14} work with the left kernel rather than the right kernel, and that in at least one of their numerical methods the point $t=0$ is ``at infinity'' with the default settings, due to a different parametrization of the kernel vectors). Setting the parameter $w = \begin{bsmallmatrix}1 & \infty & \infty & 1\end{bsmallmatrix}$ to impose the same structure as in our example, the~\texttt{slra} package of \cite{MarU14software} computes the local minimum $4$ rather than the global minimum $1$.
\end{remark}

\subsection{Smoothing the objective function}
Given $\varepsilon>0$, we now compute a relaxed version of $f(v)$ as
\begin{equation} \label{fepsilonv}
    f_\varepsilon(v) = \min_{{\Delta\in\mathcal{S}}} \,\norm{\Delta}_F^2 + \varepsilon^{-1} \norm{(A+\Delta)v}^2.
\end{equation}
\begin{theorem} \label{thm:fepsilonv}
Let a vector $v\in\mathbb{F}^n$, $v\neq 0$, and $\varepsilon>0$ be given. Define $M,r$ as in~\eqref{Mr}. Then,
\[
f_\varepsilon(v) = r^*(MM^* + \varepsilon I)^{-1}r.
\]
Moreover, the matrix that gives the minimum in~\eqref{fepsilonv} is $\Delta_* = \sum_{i=1}^p P^{(i)} (\delta_*)_i$, with
\begin{equation} \label{deltastar}
\delta_* = M^*(MM^*+\varepsilon I)^{-1}r = (M^*M+\varepsilon I)^{-1}M^*r.
\end{equation}
\end{theorem}
Note that no rank assumptions are required on $M$ here.
\begin{proof}
As in the proof of Theorem~\ref{thm:fv}, note that $\norm{\Delta}_F = \norm{\delta}$ and $(A+\Delta)v = M\delta - r$ to rewrite
\[
f_\varepsilon(v) = \min_{{\delta\in\mathbb{F}^p}} \; \norm{\delta}^2 + \varepsilon^{-1}\norm{M\delta-r}^2.
\]
Defining $\delta_*$ as in~\eqref{deltastar}, a tedious but straightforward manipulation yields the identity
\[
\delta^*\delta + \varepsilon^{-1}(M\delta -r)^*(M\delta -r) = \varepsilon^{-1}(\delta-\delta_*)^*(M^*M+\varepsilon I)(\delta-\delta_*) + r^*(MM^*+\varepsilon I)^{-1}r.
\]
Since $M^* M+\varepsilon I$ is positive definite, this formula shows that $\norm{\delta}^2 + \varepsilon^{-1}\norm{M\delta-r}^2$ is minimized when $\delta =\delta_*$, and gives the value of the minimum.
\end{proof}

We have seen in Theorem \ref{thm:fv} that the unregularized function $f(v)$ in \eqref{fv} coincides with the least-squares function $g(v)$ when $v$ is in the feasible region $\{ v: \exists \Delta \in \mathcal{S} \ s.t. \ (A+\Delta)v=0 \}$; if $A \not \in \mathcal{S}$, this feasible region can have non-empty complement in $\F^n$. On the other hand, the regularized function $f_\varepsilon(v)$ in \eqref{fepsilonv} is always defined for all $v$, and always coincides with the Tikhonov-regularized version of the least-squares function $g(v)$. In practice, we are interested in $f$ and not in $g$: When $A \not \in \mathcal{S}$, it may happen that the global minimum of $g(v)$ is strictly smaller than the global minimum of $f(v)$, but it is not achieved at a feasible $v$. Example \ref{ex:companion} below showcases an instance of this situation.

\begin{example}\label{ex:companion}
    Let $\F=\R$, and let
    \[
    {
    A=\begin{bmatrix}
        a^T\\
        e_1^T\\
        \vdots\\
        e_{n-1}^T
    \end{bmatrix}, \,
    a^T=\begin{bmatrix}
        a_{n-1} & \dots & a_1 & a_0
    \end{bmatrix}
    }
    \]
 be {a} companion matrix. We seek the singular companion matrix nearest to $A$ by allowing only perturbations of $a^T$, where the vector of perturbations is $p^T=\begin{bmatrix}
        p_{n-1} & \dots & p_1 & p_0
    \end{bmatrix}$, and without touching the ones and zeros of the companion structure. In other words, here $\mathcal{S}$ is the vector subspace of matrices of the form $\Delta=e_1 p^T$ where $p \in \R^n$, and thus clearly $A \not \in \mathcal{S}$. Of course, this example can be solved exactly without difficulty, and the unique global minimizer is $p^T=-a_0 e_n^T$; still, it is instructive to specialize our general approach to this case. The constraint $(A+ \Delta)v=0$ translates to
\[   \underbrace{e_1 v^T}_{:=M(v)} p = \underbrace{-\begin{bmatrix}
    a^T v\\
    v_1\\
    \vdots\\
    v_{n-1}
\end{bmatrix}}_{:=r(v)}. \]
Hence, the least-squares function is
\[  g(v) = \| M(v)^\dagger r(v) \|^2 = \| \overline{v} e_1^T r(v) \|^2 = |r_1(v)|^2 = |a^T v|^2 ; \]
any $v \in \spann(a)^\perp \cap \mathcal{M}$ is a minimizer of $g(v)$ and the global minimum is $0$. However, the feasible region is $F:=\spann(e_n) \cap \mathcal{M}$, and since $g(v)=|a_0|^2$ for all $v \in F$ we have that
\[  f(v)  = \begin{cases}
    |a_0|^2 \ &\mathrm{if} \ v  \in F;\\
    +\infty \ &\mathrm{if} \ v  \not\in F.
\end{cases} \]
Clearly, $\min_{{v\in\mathcal{M}}}f(v)=|a_0|^2$, coherently with our initial intuition. For the regularized function $f_\varepsilon(v)=r(v)^* (M(v)M(v)^* + \varepsilon I)^{-1} r(v)$, observe that $M(v)M(v)^* = e_1 e_1^T$, and hence,
\[  f_\varepsilon(v) = \frac{|a^T v|^2}{1+\varepsilon} + \sum_{i=1}^{n-1} \frac{|v_i|^2}{\varepsilon}.\]
It is easy to check that $f_\varepsilon(v) \rightarrow f(v)$ pointwise, when $\varepsilon \rightarrow 0$.
\end{example}

We have previously seen that, even when $A \in \mathcal{S}$, $f=g$ can be discontinuous, and therefore difficult to minimize numerically even when the feasible region is the whole $\F^n$. Instead, we sequentially minimize $f_\varepsilon(v)$, letting $\varepsilon$ tend to $0$ and using the previous minimizer as a starting point for the next iteration. The situation of Example \ref{ex:companion} is not coincidental: The sequence of functions $f_\varepsilon$ \emph{always} has the property of converging pointwise to the generalized function $f$ (but, if $A \not \in \mathcal{S}$, generally not to $g$). This is desirable in practice as we can expect that the algorithm will not be attracted towards unfeasible vectors; we state formally this result as Theorem \ref{thm:limit} below.
\begin{theorem}\label{thm:limit}
{Let the functions $f$ and $f_\varepsilon$ be defined in~\eqref{fv} and~\eqref{fepsilonv} respectively. Then, }
for all $v \in \F^n$,  $ \lim_{\varepsilon \rightarrow 0} f_\varepsilon(v) = f(v). $
\end{theorem}
\begin{proof}
   Fix $v$ and let $M(v)=\sum_{i=0}^{\min \{m,p \}} \sigma_i u_i q_i^*$ be an SVD. Then, defining $\sigma_i:=0$ when $i > \min \{m,p\}$, for all $\varepsilon > 0$ we have the eigendecomposition
   \[ (MM^* + \varepsilon I)^{-1} = \sum_{i=1}^m u_i u_i^* \frac{1}{\sigma_i^2 + \varepsilon}. \]
   Suppose first that $v$ is such that $r(v) \in \im M(v)$; then, for some vector $a$, we have $r = M a$, hence
   \[  f_\varepsilon(v) = a^* \left(\sum_{i=1}^p q_i q_i^* \frac{\sigma_i^2}{\sigma_i^2 + \varepsilon}\right)a. \]
   Since $u_i, q_i, \sigma_i, a$ do not depend on $\varepsilon$, we have
   \[ \lim_{\varepsilon \rightarrow 0} f_\varepsilon(v)= \sum_{i:\sigma_i>0} |a^* q_i|^2  =  \| M(v)^\dagger r(v)\|^2 = g(v).\]
On the other hand, if $v$ is such that $r(v) \not \in M(v)$ then for some constant (in $\varepsilon$) vectors $a$, $0 \neq b \in (\im M)^\perp$, we have $r=Ma+b$. Then, by the observations above,
\[ f_\varepsilon (v) \geq b^* (MM^* + \varepsilon I)^{-1} b = \varepsilon^{-1} \| b \|^2 \rightarrow + \infty  \ \mathrm{for} \ \varepsilon \rightarrow 0. \qedhere\]
\end{proof}

\subsection{Penalty method and augmented Lagrangian method}
\label{subsec:penalty-method-augmented-Lagrangian}

A classical algorithm to solve our minimization problem is the \emph{quadratic penalty method}~\cite[Section~4.2.1]{Bertsekas}, which consists in repeatedly solving the regularized problem, for decreasing values of $\varepsilon$. We formulate it as Algorithm~\ref{algo:penalty}.
\begin{algorithm}
\caption{Penalty method} \label{algo:penalty}
\SetKwInOut{Input}{Input}
\SetKwInOut{Output}{Output}
\Input{Initial values $\varepsilon_0, v_0 \in \mathcal{M}$; gradient tolerance $\tau$}
\Output{Approximate minimizer $v_* \approx \arg\min_{v\in\mathcal{M}} f(v)$}
$\varepsilon \gets \varepsilon_0$, $v_*\gets v_0$\;
\For{$k=1,2,\dots$}{
    $v_* \gets \arg\min_{v\in\mathcal{M}} f_\varepsilon(v)$ \tcp*{Optimization on manifolds; use the previous value of $v_*$ as starting point; stop when $\norm{\operatorname{grad} f_\varepsilon(v)} < \tau$}
    $\varepsilon \gets \varepsilon \mu$ \tcp*{$\mu<1$; either fixed or chosen adaptively}
}
\end{algorithm}

One possible issue with this method is that one obtains the exact solution only as the limit $\varepsilon \to 0$; hence, one needs to solve subproblems with smaller and smaller values of $\varepsilon$, and these subproblems become more ill-conditioned as $\varepsilon$ tends to $0$. Also, the decrease in $\varepsilon$ needs to be slow enough that at each step the previous value of $v_*$ is a suitable starting value for the subproblem with a new $\varepsilon$. For these reasons, it may be useful to consider modifications of this method for which we do not need to push the parameter $\varepsilon$ to $0$.

Potential help can come from the so-called \emph{augmented Lagrangian method}~\cite{Bertsekas}, also known as \emph{method of multipliers}. Combining Riemannian optimization with augmented Lagrangian methods for some specific matrix nearness problems was proposed, in a different setting than in this paper, in \cite{borsdorf}.

In an augmented Lagrangian method, the main idea is combining the penalty method with a Lagrangian formulation. To introduce the method, we start from the constrained minimum problem on the pair $(\Delta, v)$
\[
\min_{{(\Delta,v) \in \mathcal{S} \times \mathcal{M}}} \norm{\Delta}_F^2 \subjectto (A+\Delta)v = 0.
\]
We have shown in Theorem~\ref{thm:fv} that we can obtain  the optimal $\Delta$ in closed form for any given $v$, hence reducing the minimization to $\min_v f(v)$. Instead of doing that, we first switch to the augmented Lagrangian formulation of this problem, following~\cite[Section~4.2]{Bertsekas}: We introduce a dual variable $y\in\mathbb{F}^{m}$, and consider the Lagrangian function
\[
\mathcal{L}_\varepsilon(\Delta,v;y) = \norm{\Delta}_F^2 + \varepsilon^{-1}\norm{(A+\Delta)v}^2 + 2\left\langle y, (A+\Delta)v \right\rangle.
\]
In the augmented Lagrangian method, at each step $k$ we compute (for fixed $\varepsilon$ and $y$) the minimum
\[
\min_{\Delta,v} \mathcal{L}_{\varepsilon}(\Delta,v;y).
\]
Even in this generalized setting, for every given $v$, the minimum over $\Delta$ can still be obtained in closed form. To this effect, note that (since $y$ is fixed) this new minimization task is equivalent to
\begin{align*}
    \min_{\Delta,v} \mathcal{L}_{\varepsilon}(\Delta,v;y) + \varepsilon \norm{y}^2 &= {\min_{\Delta,v}} \; \norm{\Delta}_F^2 + \varepsilon^{-1}\left\langle  (A+\Delta)v+\varepsilon y, (A+\Delta)v+\varepsilon y\right\rangle\\
    &= {\min_{\Delta,v}} \; \norm{\Delta}_F^2 + \varepsilon^{-1}\norm{(A+\Delta)v+\varepsilon y}^2.
\end{align*}
This function is very similar to $f_\varepsilon(v)$ in~\eqref{fepsilonv}, with the only difference that we have replaced $Av$ with $Av + \varepsilon y$. Hence, we can replicate the proof of Theorem~\ref{thm:fepsilonv} almost verbatim, with the only difference being the substitution of $r=-Av$ with $r=-Av-\varepsilon y$. This way, we obtain
\begin{equation*}
    \begin{gathered}
        f_{\varepsilon, y}(v) = \min_{\Delta\in\mathcal{S}} \mathcal{L}_{\varepsilon}(\Delta,v;y) + \varepsilon \norm{y}^2 = r^*(MM^*+\varepsilon I)^{-1}r,\\ M=M(v),\quad r = -Av - \varepsilon y. 
    \end{gathered}
\end{equation*}
The resulting method is summarized in Algorithm~\ref{algo:auglag}.
\begin{algorithm}
\caption{Augmented Lagrangian method} \label{algo:auglag}
\SetKwInOut{Input}{Input}
\SetKwInOut{Output}{Output}
\Input{Initial values $\varepsilon_0, v_0 \in \mathcal{M}, y_0\in\mathbb{C}^m$; gradient tolerance $\tau$}
\Output{Approximate minimizer $v_* \approx \arg\min_{v\in\mathcal{M}} f(v)$}
$\varepsilon \gets \varepsilon_0$, $v_*\gets v_0$\;
\For{$k=1,2,\dots$}{
    $v_* \gets \arg\min_{v\in\mathcal{M}} f_{\varepsilon,y}(v)$ \tcp*{as in Algorithm~\ref{algo:penalty}}
    $\Delta_* \gets \Delta_*(v_*)$ \tcp*{value of the minimizer}
    $y \gets y + \varepsilon^{-1} {(A+\Delta_*)v_*}$ \tcp*{update $y$ with gradient ascent}
    $\varepsilon \gets \varepsilon \mu$ \tcp*{as in Algorithm~\ref{algo:penalty}}
}
\end{algorithm}
Algorithm \ref{algo:auglag} has the potential for faster convergence, since there are two factors that drive the iterates close to the optimal solution: The convergence of $\varepsilon$ to $0$, and the convergence of $y$ to the optimal value of the Lagrange multiplier $y_*$. We refer the reader to the discussion in \cite[Section~4.2.2]{Bertsekas} for more details. In practice, in our numerical experiments, we have found that Algorithm \ref{algo:penalty} can sometimes be competitive when cleverly implemented. Moreover, we experienced that employing a suitable update strategy for the regularization parameter may avoid the need to solve ill-conditioned subproblems, for small values of $\varepsilon$; see Subsection \ref{subsec:adaptive strategy} for the implementation details.

In both the penalty method and the augmented Lagrangian approach, we solve minimization problems on Riemannian manifolds, producing a sequence of minimizers. It is useful to discuss convergence properties of this sequence; Section \ref{subsec:convergence-results} provides a theoretical convergence analysis.

We note that the software in~\cite{MarU14software} also uses a form of regularization, but a different, unrelated one: they relax the constraint $v^*v = 1$ (or more generally $V^*V=I$ in higher-nullity problems) as a substitute for Riemannian optimization; this is not the same as relaxing the constraint $(A+\Delta)v = 0$ as we do here.
Another method to minimize over the Grassmann manifold is suggested in~\cite{UseM14manifold}, but once again without the regularization described here.

It remains to discuss how to solve the inner subproblems $\min_v f_{\varepsilon,y}(v)$; we do this in the next section.

\section{Minimization of $f(v)$ via Riemannian optimization}\label{sec:riemann}

In the algorithms introduced in the previous section, each inner subproblem requires minimizing $f_{\varepsilon,y}(v)$ (for some $\varepsilon,y$; the penalty method in Algorithm \ref{algo:penalty} can be seen as a special case of Algorithm \ref{algo:auglag} with $y=0$) over the unit sphere $\mathcal{M}$, which is a real embedded Riemannian manifold (even in the case $\F$ = $\C$). To this end, we may employ the MATLAB package Manopt \cite{BoumalMishraAbsil}, a toolbox to solve optimization problems on manifolds. Specifically, we use its Riemannian trust-region method; see~\cite{AbsilBaker} for details. Manopt's trust-region method requires providing the Riemannian gradient and the Riemannian Hessian of the function $f_{\varepsilon,y}\big|_{\mathcal{M}}: \mathcal{M} \mapsto \mathbb{R}$. 

The Riemannian gradient can be computed as the projection onto the tangent space $T_v(\mathcal{M})$ of the Euclidean gradient of $f_{\varepsilon,y}$ in the ambient space $\mathbb{F}^n$; see \cite[Section~3]{Boumal}.
The Riemannian Hessian, instead, can be computed from the Euclidean gradient and Hessian through the Weingarten map for the manifold $\mathcal{M}$ and the orthogonal projection onto the tangent space $T_v (\mathcal{M})$; see \cite[Section 5]{Boumal}. In the special case of the unit sphere, explicit formulae are given in \cite{AbsilMahony}. The construction of the Riemannian gradient and the Riemannian Hessian is handled automatically in Manopt, which provides functions for its computation from their Euclidean counterparts. Therefore, we only need to provide expressions for the Euclidean gradient and the Euclidean Hessian-vector product. In~\cite[Section~4]{UseM14}, the authors obtain (using different notation) an expression for the Euclidean gradient of $f(v)$ and an approximation of its Hessian; here, we improve on that by expanding the computation to $f_{\varepsilon,y}(v)$ and giving the exact Hessian.

\begin{theorem} \label{thm:gradienthessian}
    Following the above notation, let $M = M(v)$ and
    \begin{align*}
        r &= -Av - \varepsilon y &  z &= (MM^*+\varepsilon I)^{-1}r, \\\delta_* &= M^*z, & \Delta_* &= \sum_{i=1}^p P^{(i)} (\delta_*)_i.
    \end{align*}
    Then, the Euclidean gradient of $f_{\varepsilon,y}(v) = r^*(MM^*+\varepsilon I)^{-1}r$ is
    \begin{equation} \label{gradf}
         \nabla f_{\varepsilon,y}(v) = -2(A+\Delta_*)^*z,
    \end{equation}
    and the Hessian $H_{\varepsilon,y}(v)$ is the matrix whose action is
    \[
        H_{\varepsilon,y}(v) w = -2 (\dot{\Delta}_*)^*z - 2(A+\Delta_*)^*\dot{z},
    \]
    with
    \begin{gather}
    \dot{z} = -(MM^*+\varepsilon I)^{-1}(M\dot{M}^*z + (A+\Delta_*)w), \label{dotz}\\
        \begin{aligned}
        \dot{M} &= M(w), &
        \dot{\delta}_* &= \dot{M}^* z + M^*\dot{z}, &
        \dot{\Delta}_* &= \sum_{i=1}^p P^{(i)} (\dot{\delta}_*)_i. \label{dotstuff}
        \end{aligned}
    \end{gather}
\end{theorem}

The proof of Theorem \ref{thm:gradienthessian} is in Appendix~\ref{sec:proofgradienthessian}.

\subsection{Computational cost}

Let us briefly indulge in some observations on the computational complexity of our algorithm. The most expensive operations needed to compute $f_{\varepsilon,y}(v), \nabla f_{\varepsilon,y}(v)$ and $H_{\varepsilon,y}(v)w$ are:
\begin{enumerate}
    \item A matrix-vector product with $A$ and one with $A^*$, for a cost of $\mathcal{O}(mn)$ (in the most general case; this cost can be reduced by exploiting matrix structures such as sparsity or block Toeplitz/Hankel, when present).
    \item Solving two linear systems with the matrix $MM^* + \varepsilon I$, for a cost of $\mathcal{O}(mp \min(m,p))$ if one does it using an SVD of $M$; this is beneficial to stability anyway, since we have seen that ill-conditioning and singularity of $M$ are a major concern.
    \item Computing (twice each) $M(x)$ and $P(y) = \sum_{i=1}^p P^{(i)}y_i$ for given vectors $x\in\mathbb{F}^n$ and $y\in\mathbb{F}^{p}$; these operations cost $\mathcal{O}(mnp)$ (again, ignoring any structure in the $P^{(i)}$ that could be exploited to reduce the cost).
\end{enumerate}
The total cost $\mathcal{O}(mp(n+\min(m,p)))$ must be then multiplied by the (possibly large) number of evaluations needed in Algorithm~\ref{algo:penalty} or Algorithm~\ref{algo:auglag}. In the next sections, we will see that these opportunities for reduced complexity do in fact arise in several practically relevant cases.

\section{Convergence results}
\label{subsec:convergence-results}

As described in Subsection \ref{subsec:penalty-method-augmented-Lagrangian}, and still noting that the penalty method is a special case of the more general augmented Lagrangian method, our approach relies on the solution of the unconstrained minimization problem
\[
\min_{(\delta, v) \in \mathbb{F}^p \times \mathcal{M}} \| \delta \|^2 + \varepsilon^{-1}\| Av + M(v) \delta\|^2 + 2 \left\langle y, Av + M(v) \delta \right\rangle = \min_{(\delta, v) \in \mathbb{F}^p \times \mathcal{M}}\mathcal{L}_\varepsilon(\delta,v;y) 
\]
to provide a numerical approximation of the solution of the constrained minimization problem
\begin{equation*}
    \min_{{(\delta, v) \in \mathbb{F}^p \times \mathcal{M}}} \| \delta \|^2 
     \subjectto  Av + M(v) \delta =0.
\end{equation*}
Consider a sequence of parameters ($\varepsilon_k$, $y_k$) for the augmented Lagrangian formulation and denote by
\begin{equation*}
\label{eq:step_augmented_lagrangian}
(\delta_k, v_k) = \argmin_{(\delta, v) \in \mathbb{F}^p \times \mathcal{M}} 
 \mathcal{L}_{\varepsilon_k}(\delta,v;y_k)
\end{equation*}
the solution of this subproblem. We follow the methodology in \cite{LiuBoumal} and apply  \cite[Proposition 3.2]{LiuBoumal} to the sequence of solutions $( \delta_k, v_k)$. For completeness, we state the convergence result below as Proposition \ref{prop:convergence}. For the technical definitions of Linear Independence Constraint Qualifications (LICQ) and First-Order Necessary Conditions (KKT conditions) in the setting of Riemannian manifolds, we refer to~\cite[Equations (4.3) and (4.8)]{Zhang}.
\begin{proposition}[\cite{LiuBoumal}] \label{prop:convergence}
Let $\Omega$ denote the set of feasible points. Applying the Riemannian augmented Lagrangian method  \cite[Algorithm 1]{LiuBoumal} with a sequence $\eta_k \rightarrow 0$, if at each iteration $k$ the subsolver produces a point $x_{k+1}$ satisfying
\begin{equation}
\label{eq:stop_gradient}
    \left\|\operatorname{grad}_x \mathcal{L}_{\varepsilon_k}\left(x_{k+1}, y^k\right)\right\| \leq \eta_k,
\end{equation}
and if the sequence $\left\{x_k\right\}_{k=0}^{\infty}$ has a limit point $\bar{x} \in \Omega$ where LICQ conditions are satisfied, then $\bar{x}$ satisfies KKT conditions of the original constrained minimization problem. 
\end{proposition}

Proposition \ref{prop:convergence} states, in the general case of Riemann {manifolds}, that under the assumption of LICQ the limit point satisfies the KKT conditions, which is the typical desired result in optimization. In the more classical context of optimization on Euclidean spaces, analogous convergence results are better known; see for instance \cite{Birgin} for a complete overview.

At each step of our approach, we solve an unconstrained minimization problem in the form \eqref{eq:step_augmented_lagrangian}, where in our setting, the ordered pair $(\delta, v)$ takes the place of the variable $x$ in Proposition \ref{prop:convergence}. Our method can be seen as a partially explicit method, since we have an explicit formula for  the optimal value of $\delta$, and a partial trust-region method on the Riemannian manifold $\mathcal{M}$ for the variable $v$. In our problem, the LICQ conditions reduce to assuming that the matrix
\[
\begin{bmatrix}
    M(v)^* \\
    (I - vv^*) (A + \Delta)^*
\end{bmatrix} \in \mathbb{F}^{(n+p) \times n}
\]
has full rank at a limit point $(\bar{\delta}, \bar{v})$ of the sequence $\left\lbrace ( \delta_k, v_k) \right\rbrace_k$. Then, the theorem shows that $(\bar{\delta}, \bar{v})$ satisfies the KKT conditions for the constrained minimization problem. More convergence results exist for the employment of second-order methods as inner solvers at each iteration $k$ (\cite[Proposition $3.4$]{LiuBoumal}), as long as all the eigenvalues of  $\operatorname{Hess}_x \mathcal{L}_{\varepsilon_k}\left(x_{k+1}, y^k\right)$ are equal to or above the value $-\eta_k$. Note that, in our implementation of the method, we propose a stopping criterion only based on the condition \eqref{eq:stop_gradient}, and therefore we can only guarantee convergence towards stationary points.

The same convergence analysis applies to the penalty method, starting with the Lagrangian multiplier $y=0$ and not updating it (see for instance the update proposed in \cite[Algorithm $1$]{LiuBoumal}). Additionally, we observe that convergence results for the Riemannian trust-region method used in the inner problems can be found in \cite[Section $4$]{AbsilBaker}.

\section{Nearest sparse singular matrix} \label{sec:sparse}
In this section, we demonstrate how the computations can become simpler if the specific structure of the vector space $\mathcal{S}$ is exploited. We illustrate this by a common case, where $\mathcal{S}$ is the set of matrices with a prescribed sparsity pattern. Given a set $\mathcal{J} \subseteq \{1,2,\dots,n\}^2$, we are interested in perturbations $\Delta$ such that $\Delta_{ij}=0$ whenever $(i,j) \not \in \mathcal{J}$. A natural basis for $\mathcal{S}$ is
\[
\{P^{(1)}, \dots, P^{(\# \mathcal{J} )}\} = \{e_i e_j^* \colon (i,j) \in \mathcal{J}\}.
\]
In particular, the columns of $M$ are of the form $e_i v_j$ for all $(i,j)\in \mathcal{J}$ (taken in a prescribed order), and $MM^*$ is the diagonal matrix
\[
MM^* = \operatorname{diag}(d_1,\dots, d_n), \quad d_{i} = \sum_{(i,j)\in \mathcal{J}} \abs{v_j}^2,\text{ for all $i=1,\dots,n$}.
\]
This leads to expressions that can be used to compute $f_{\varepsilon, y}$, its gradient and the action of its Hessian in $\mathcal{O}(\# \mathcal{J})$ operations; for instance
\[
    f_{\varepsilon,y}(v) = \norm{(Av+\varepsilon y) \odot h}^2, \quad h_i = \frac{1}{\sqrt{d_i^2+\varepsilon}},
\]
\[
\nabla f_{\varepsilon,y}(v) = 2A^*(h \odot (Av+\varepsilon y)) - 2 (f \odot v), \quad f_j = \sum_{(i,j)\in \mathcal{J}} \abs{z_i}^2,
\]
where $\odot$ denotes the Hadamard (componentwise) product of vectors.

When $p\geq n$, we can also find explicitly a thin SVD of $M$ as $M = I \cdot S \cdot (S^{-1}M)$, with $S = \operatorname{diag}(d_1^{1/2}, \dots, d_n^{1/2})$.

\section{Nearest singular matrix polynomial}
\label{sec:Nearest_sing_matrix_poly}

In some applications, it is of interest to find the distance to singularity for a given matrix polynomial $A(x) \in \mathbb{F}[x]^{n\times n}_k$, and to find the corresponding nearest singular matrix polynomial \cite{nearsingpen, bora, DNN24}. 
The distance between two matrix polynomials is usually measured using the distance induced by the norm 
\[
\norm{ A(x)}_F:=\norm*{\begin{bmatrix} A_0 & A_1 & \ldots  & A_k \end{bmatrix}}_F.
\]

The coefficients of the matrix polynomial usually arise from a model of a physical system, and in many cases, such as in \cite{Grabner}, the grade (maximal degree) of the matrix polynomial is fixed. Therefore, in this section we focus on finding the closest singular matrix polynomial with the same grade as that of $A(x)$. However, it is worth noting that this does not always coincide with the nearest singular matrix polynomial of unbounded degree, as shown below in Example \ref{ex:nearest_poly}.

\begin{example} \label{ex:nearest_poly}
Let $0 < \varepsilon < 1$ and
    \[ P(x) = \begin{bmatrix}
        0 & x\\
        x & \frac{1}{\varepsilon}
    \end{bmatrix}.\]
The distance of $P(x)$ from the set of singular pencils (polynomials of degree at most $1$) can be computed exactly and is equal to $1$. However, the degree $2$ polynomial
\[ Q(x) = \begin{bmatrix}
    \varepsilon x^2 & x\\
    x & \frac{1}{\varepsilon}
\end{bmatrix} \]
is singular and only $\varepsilon < 1$ away from $P(x)$. Moreover, this same example also shows that the ratio of the two distances (with constrained or unconstrained degree) is unbounded.    
\end{example}

Let us now formalize the problem. Given a matrix polynomial $A(x) \in \mathbb{F}[x]^{n\times n}_k$ in the form $A(x)=A_0 + A_1 x + \ldots + A_k x^k$, where $A_i \in \mathbb{F}^{n \times n}$ for $i=0,\ldots,k$, we want to find the smallest perturbation $\Delta(x) \in \mathbb{F}[x]^{n\times n}_k$ that gives a singular matrix polynomial.
A matrix polynomial $A(x) \in \mathbb{F}[x]^{n\times n}_k$ is called \emph{singular} if any of the following  equivalent conditions holds:
\begin{itemize}
\item $\det A(x) = 0$;
    \item There exists a nonzero rational vector $v(x) \in \mathbb{F}(x)^{n}$ such that $A(x)v(x) = 0$;
    \item There exists a nonzero vector polynomial $v(x)=v_0 + v_1 x + \ldots + v_d x^d$, with $d \leq k(n-1)$, such that $A(x)v(x) = 0$;
    \item There exists a nonzero vector polynomial $v(x)$, with $d \leq \lfloor \frac{k(n-1)}{2}\rfloor$, such that either $A(x)v(x)=0$ or $v(x)^* A(x) =0$.
\end{itemize}
The last {two conditions, which are less obvious, are} a consequence of Theorem \ref{thm:orders}, whose statement and discussion is postponed to a later section.

For a fixed $d$, we may define
\begin{equation}
\label{eq:Toeplitz matrix_for matrix polynomials}
    \mathcal{T}_d(A) = \begin{bmatrix}
    A_0\\
    A_1 & A_0\\
    \vdots & A_1 & \ddots\\
    A_k & \vdots & \ddots & A_0\\
    & A_k & \ddots & A_1\\
    & & \ddots & \vdots\\
    & & & A_k
\end{bmatrix} \in \mathbb{F}^{n(k+d+1) \times n(d+1)},
\quad
\vvec v = \begin{bmatrix}
    v_0\\
    v_1\\
    \vdots\\
    v_d
\end{bmatrix} \in \mathbb{F}^{n(d+1)}.
\end{equation}
Note that in general $\vvec \left( A(x) v(x) \right)= \mathcal{T}_d(A) \vvec v$. {Indeed,} $A(x)v(x)=\sum_{i=0}^{d+k} \left( \sum_{j+l=i} A_j v_l \right) x^i$. Since we are interested in looking for a singular matrix polynomial, we use the relation 
\[
A(x) v(x) = 0 \iff \mathcal{T}_d(A) \vvec v = 0.
\]
Hence we can convert the problem of finding the nearest singular matrix polynomial to $A(x)$ into finding the nearest singular matrix to $\mathcal{T}_d(A)$ that has the same block Toeplitz structure.

Linear perturbation structures can be converted analogously. Given matrix polynomials { $P^{(1)}(x)$, $ \dots, P^{(p)}(x)$},
\[
\Delta(x) \in \mathcal{S} = \left\{\sum_{i=1}^p P^{(i)}(x) \delta_i \colon \delta\in\mathbb{F}^p\right\}
\]
holds if and only if
\[
\mathcal{T}_d(\Delta) \in \left\{\sum_{i=1}^p \mathcal{T}_d(P^{(i)}) \delta_i \colon \delta\in\mathbb{F}^p \right\} =: \mathcal{T}_d(\mathcal{S}).
\]
The norm is also preserved up to a constant:
\[
\norm{\mathcal{T}_d(\Delta)}_F = \sqrt{d+1} \norm*{\Delta(x)}_F.
\]

Then the problem of finding the nearest singular matrix polynomial to $A(x)$ of degree $\leq k$, with a {nonzero} polynomial $v(x)$ of degree at most $d$ in its kernel 
\begin{equation} \label{eq:unregular_polynomial_problem}
    \widetilde{f}(v):= \min_{{\Delta (x) \in \mathcal{S} }} 
    \norm{\Delta(x)}^2_F \subjectto \left( A(x) + \Delta(x) \right) v(x) =0, 
\end{equation}
can be rephrased into finding the nearest singular structured matrix to $\mathcal{T}_d(A)$ with $\vvec v$ in its kernel:
\[
f(v):=\min_{{\mathcal{T}_d(\Delta)  \in \mathcal{T}_d(\mathcal{S})}} 
\norm{\mathcal{T}_d(\Delta)}^2_F \subjectto \left( \mathcal{T}_d(A) + \mathcal{T}_d(\Delta) \right)\vvec v =0  
\]
The minimizers for $f(v)$ and $\widetilde{f}(v)$ are the same (up to representation), while the mimima differ by a factor $d+1$.

In this form, it is clear that the problem is a special case of the (structured) nearest singular matrix problem given in \eqref{fv}: The matrix under consideration is $\mathcal{T}_d(A)$ and the set of structured perturbations is $\mathcal{T}_d(\mathcal{S})$. In order to find the distance to singularity, we minimize $f(v)$ with the constraint $\vvec v \in \mathcal{M} \subseteq \mathbb{F}^{n(d+1)}$. With these slight notation changes, we can tackle the polynomial variant of the problem with the same machinery that was developed in Section \ref{sec:matrix} for the matrix case.

\begin{example}
We consider a quadratic matrix polynomial $A(x) \in \mathbb{R}[x]_2^{2\times 2}$ in the form
\[
A(x)=x^2 \begin{bmatrix}
    1  & 0 \\
    0 & 0 
\end{bmatrix} + x \begin{bmatrix}
    0 & 1 \\
    2 & 0
\end{bmatrix} + \begin{bmatrix}
    0 & 0 \\
    0 & 1\\
\end{bmatrix},
\]
and we impose the sparsity pattern induced by the coefficient matrices. With the previous notations, we have:
\[
P^{\left(1\right)}(x)= x^2 \begin{bmatrix}
    1& 0\\
    0 & 0
\end{bmatrix}, P^{\left(2\right)}(x)= x \begin{bmatrix}
   0 & 1\\
    0 & 0
\end{bmatrix}, P^{\left(3\right)}(x)= x \begin{bmatrix}
    0& 0\\
    1 & 0
\end{bmatrix}, P^{\left(4\right)}(x)= \begin{bmatrix}
 0& 0\\
    0 & 1
\end{bmatrix}.
\]
Setting $d=2$ and $p=4$, the set $\mathcal{T}_2(\mathcal{S})$ is given by the following basis:
\[
\mathcal{T}_2(P^{\left(1\right)}) = \begin{bmatrix}
    0 \\
    0 &0 \\
    \begin{bmatrix}
        1 & 0 \\
        0 & 0
    \end{bmatrix} & 0 & 0\\
     &  \begin{bmatrix}
        1 & 0 \\
        0 & 0
    \end{bmatrix} & 0 \\
    & & \begin{bmatrix}
        1 & 0 \\
        0 & 0
    \end{bmatrix} 
\end{bmatrix}, \; \mathcal{T}_2(P^{\left(2\right)}) = \begin{bmatrix}
    0 \\
     \begin{bmatrix}
        0 & 1\\
        0 & 0
    \end{bmatrix}  &0 \\
   0 &  \begin{bmatrix}
        0 & 1 \\
        0 & 0
    \end{bmatrix} & 0\\
     &  0 & \begin{bmatrix}
        0 & 1 \\
        0 & 0
    \end{bmatrix} \\
    & & 0 
\end{bmatrix},
\]
\[
\mathcal{T}_2(P^{\left(3\right)}) = \begin{bmatrix}
    0 \\
     \begin{bmatrix}
        0 & 0\\
        1 & 0
    \end{bmatrix}  &0 \\
   0 &  \begin{bmatrix}
        0 & 0 \\
        1 & 0
    \end{bmatrix} & 0\\
     &  0 & \begin{bmatrix}
        0 & 0 \\
        1 & 0
    \end{bmatrix} \\
    & & 0 
\end{bmatrix}, \; \mathcal{T}_2(P^{\left(4\right)}) = \begin{bmatrix}
    \begin{bmatrix}
        0 & 0\\
        0 & 1
    \end{bmatrix} \\
     0  &  \begin{bmatrix}
        0 & 0\\
        0 & 1
    \end{bmatrix} \\
   0 & 0 &  \begin{bmatrix}
        0 & 0\\
        0 & 1
    \end{bmatrix}\\
     &  0 & 0\\[1mm]
    & & 0 
\end{bmatrix}.
\]
\end{example}

\subsection{Unstructured matrix polynomials}
\label{subsec:unstructured_distance_polynomial}

The framework becomes simpler if we do not impose additional structures on the matrix polynomial $\Delta(x)$. This leads to a computationally more efficient algorithm that will be subsequently used for numerical experiments in Subsection \ref{sec:betterthanbora}. In this {section}, we explain the formulation of this more efficient algorithm that uses the penalty method approach, presented in Subsection \ref{subsec:penalty-method-augmented-Lagrangian}.

The relation $(\mathcal{T}_d(A) + \mathcal{T}_d(\Delta)) \vvec v=0$ can be reformulated as
\[
\begin{bmatrix}
\Delta_0 v_0 \\
\Delta_1 v_0 + \Delta_0 v_1 \\
\vdots \\
\sum_{j+l=i}\Delta_j v_l \\
\vdots \\
\Delta_k v_d
\end{bmatrix}=-
\begin{bmatrix}
A_0v_0 \\
A_1 v_0 + A_0 v_1 \\
\vdots \\
\sum_{j+l=i}A_j v_l \\
\vdots \\
A_k v_d
\end{bmatrix},
\]
or equivalently
\begin{equation} \label{eq:unstr_matrix_pol_relation}
\begin{bmatrix}
    \Delta_0 & \Delta_1 & \ldots & \Delta_k
\end{bmatrix} \mathcal{T}_{k}(v^T)^T = - \begin{bmatrix}
    A_0 & A_1 & \ldots & A_k
\end{bmatrix} \mathcal{T}_{k}(v^T)^T,
\end{equation}
where $\mathcal{T}_{k}(v^T) \in \mathbb{F}^{(k+d+1) \times n(k+1)}$ is defined as in \eqref{eq:Toeplitz matrix_for matrix polynomials}.
Using the properties of the Kronecker product, the relation \eqref{eq:unstr_matrix_pol_relation} can be written as
\[
\left( \mathcal{T}_k(v^T) \otimes I_n \right)\delta = - \left( \mathcal{T}_k(v^T) \otimes I_n \right) \alpha,
\]
where $\delta:=\vvec{\left(\begin{bmatrix}
    \Delta_0 & \ldots & \Delta_k
\end{bmatrix}\right)}$ and $\alpha:=\vvec{\left(\begin{bmatrix}
    A_0 & \ldots & A_k
\end{bmatrix}\right)}$.
Hence, we can set
\begin{equation} \label{Mrpolynomial}
    M(v)=\left( \mathcal{T}_k(v^T) \otimes I_n\right), \quad r(v)=-M(v) \alpha,
\end{equation}
and apply the results in the previous {sections}. In particular, Theorem~\ref{thm:fepsilonv} gives
\begin{align*}
f_{\varepsilon}(v) &= \alpha^*M(v)^* \left( M(v) M(v)^* + \varepsilon I_{n(k+d+1)}\right)^{-1} M(v) \alpha,
\end{align*}
and after some tedious manipulation, this can be expressed as
\begin{align*}
f_{\varepsilon}(v) &= \left\langle \begin{bmatrix}
    A_0 & \ldots & A_k
\end{bmatrix} , \begin{bmatrix}
    A_0 & \ldots & A_k
\end{bmatrix} \mathcal{T}_k(v^T)^T  C(v)^T \mathcal{T}_k(v^*)  \right\rangle,
\end{align*}
where $C(v):=(\mathcal{T}_k(v^T) \mathcal{T}_k(v^T)^* +\varepsilon I_{(k+d+1)} )^{-1}$.

Similarly, one can use the Kronecker product structure in $M(v)$ from~\eqref{Mrpolynomial} to reduce the computational cost in the formulae for the (Euclidean) gradient and the (Euclidean) Hessian in Theorem~\ref{thm:gradienthessian}. Indeed, the square matrix $\mathcal{T}_k(v^T)^T C(v)^T\mathcal{T}_k(v^*)$ has size $n(k+1)$, while the matrix $M(v)$ is $n(k+d+1) \times n^2(k+1)$.

As observed in \cite{bora}, we can reduce the sizes of the matrices even further by searching through the left and the right kernel separately, with $d = \lfloor \frac{k(n-1)}{2}\rfloor$, instead of minimizing \eqref{eq:unregular_polynomial_problem} with $d = k(n-1)$. In this way, we need to run the algorithm twice, but the sizes of the matrices are significantly smaller. Overall, this change  translates into faster computational time.  

\subsection{Discontinuities}
In Subsection \ref{sec:need_regularization}, we saw that $f$ can become discontinuous if structure is imposed on the problem. As discussed earlier in this section, the nearest singular matrix polynomial problem introduces a layer of structure given by the operator $\mathcal{T}_d$. It turns out that this causes discontinuities to emerge even when no additional structure is imposed on the coefficients of $\Delta(x)$. Moreover, we can gain insight into the locations of these discontinuities by applying a matrix-theoretical framework. In particular, we invoke the concept of minimal indices. A complete treatment of minimal indices is a somewhat technical topic in the theory of matrix polynomials, and beyond the scope of our paper. Interested readers can consult for example \cite{dd,forney,noferini24} and the references therein. Here we only treat the special case of a singular $n \times n$ matrix polynomial $A(x)$ of rank $n-1$. In this  case, it is clear that there exists a polynomial vector $v(x)$ such that the (free) module $\ker A(x) \cap \F[x]^n = \mathrm{span}_{\F[x]} v(x)$. It turns out that $\deg v(x)$ is uniquely determined by $A(x)$ and it is called the right minimal index of $A(x)$; a left minimal index is defined analogously.

The following Theorem \ref{thm:orders} is important for our next developments; we omit its proof because it is a corollary of several more general results including the Index Sum Theorem; {see} \cite{ist} and \cite[Theorem 3.2]{dd}.

\begin{theorem}\label{thm:orders}
    Let $A(x) \in \F[x]^{m \times n}$ be a polynomial matrix having degree ${\leq} d$ and rank ${\leq} n-1$. Denote the {sum of the} left and right minimal indices of $A(x)$ by $\eta_\ell$ and $\eta_r$ respectively. Then, $\eta_\ell + \eta_r \leq d{(n-1)}$. Moreover, equality holds generically in the set of $n \times n$ polynomial matrices of degree $\leq d$ and rank ${\leq} n-1$.
\end{theorem}

\begin{corollary}\label{cor:nowheredense}
    Let $A(x) \in \F[x]^{n \times n}$ be a singular square polynomial matrix of degree $d$. Then, generically {in the sense of Theorem \ref{thm:orders}}, $A(x)$ has rank $n-1$, right minimal index $\eta_r$ and left minimal index $\eta_\ell$ satisfying $\eta_r+\eta_\ell = d(n-1)$.
\end{corollary}

\begin{proof}
Denote by $S$ the set of $n \times n$ polynomial matrices of degree $\leq d$.  For  $k=0,1,\dots,n-1$, let $S_k \subsetneq S$ be the subset of $n \times n$ polynomial matrices of degree $\leq d$ and rank $\leq k$. By the definition of rank,
    \[ S_{n-1}=\{ A(x) \in S : \det A(x) = 0 \}, \quad   S_{n-2}=\{ A(x) \in S : \mathrm{adj}A(x)=0 \} \subsetneq S_{n-1}. \]
We only need to establish that $S_{n-2}$ is nowhere dense in $S_{n-1}$, as the rest of the statement then follows immediately by using Theorem \ref{thm:orders} {and the fact that the intersection two open dense sets is open and dense}. It is clear that $S_{n-2}$ is closed in $S_{n-1}$, by continuity of determinants. Thus, it suffices to show that, for all $A(x) \in S_{n-2}$ and for all $\delta >0$, the ball of radius $\delta$ and centered at $A(x)$ contains at least one element $B(x) \in S_{n-1}\setminus S_{n-2}$. Let $\lambda \in \F$ be such that $\rank A(\lambda) = \rank A(x)=r < n-1$. Then by standard results in matrix theory (for example using the SVD), one can find $U,V \in \F^{n \times (n-1-r)}$ such that $\|UV^T\|_F = 1$ and $A(\lambda) + \delta U V^T$ has rank $n-1$. On the other hand,
    \[ \rank \left( A(\lambda) + \delta U V^T \right) \leq \rank \left(A(x) + \delta U V^T\right) \leq \rank\left(A(x)\right) + \rank \left(\delta U V^T \right)  \]
    and hence
    \[ n-1 \leq \rank \left(A(x) + \delta U V^T\right) \leq r + n-1-r=n-1.  \]
\end{proof}

It is also clear\footnote{For example if $A(x)=\sum_{i=0}^d A_i x^i$ ($A_d=0)$, we can set
\[[A(x)]_k = x^d \sum_{i=0}^d A_{d-i} \left(\frac{1}{x}+\frac{a_k}{k}\right)^i    \]
 where $(a_k)_k$ is any sequence such that $0<|a_k| \leq 1$ and $0 \neq A\left(\frac{k}{a_k}\right)$ for all $k$.} that, given a singular $n \times n$ polynomial $A(x) \neq 0$ of degree $\leq d$, there is a sequence $A_k(x)$ of singular $n \times n$ polynomials of degree exactly $d$ s.t. $A_k(x) \rightarrow A(x)$. {An immediate consequence of this fact and of Corollary \ref{cor:nowheredense} is that, given a regular $n \times n$ $R(x)$ of degree $d$ and a singular polynomial $R(x)+E(x)$ nearest to $R(x)$ and of degree $\leq d$, then there exist another singular polynomial $R(x)+F(x)$ of the same degree $d$, whose rank is precisely $n-1$, whose left and right minimal indices sum to $d(n-1)$, and such that $F(x)$ is arbitrarily close to $E(x)$.}

\begin{lemma}\label{lemma:rdrop_minindex}
      The matrix $\mathcal{T}_{k}(v^T)$ as in \eqref{eq:unstr_matrix_pol_relation} has left nullity $\geq d-m$ if $v(z)=v_0+v_1 z + \dots + v_{d} z^d$ has the form $v(z)=p(z) w(z)$ where $w(z)$ is a minimal basis of degree $m \leq d$ and $p(z)$ is a scalar polynomial of degree $\leq d-m$.
\end{lemma}
  
\begin{proof}
Note that $\mathcal{T}_{k}(v^T)$ has $k+d+1$ rows and $(k+1)n$ columns.

    Write $p(z) = \prod_{j=1}^{d-m} (a_j z - b_j)$, where $a_j,b_j$ are not both $0$ (if, for some $j$, $a_j=0$ then the degree of $p(z)$ is strictly less than $d-m$, and we see it as having some roots at infinity). Then, $p(z)$ has $d-m$ homogeneous roots (including possibly infinity) at the points of homogeneous coordinates $(b_j,a_j)$. Note that homogeneous roots can be repeated according to their multiplicity.

Observe that a vector of the form

\[  w(a,b)= \begin{bmatrix}
    a^{k+d} & a^{k+d-1} b & a^{k+d-2} b^2 & \dots & a b^{k+d-1} & b^{k+d}
\end{bmatrix}    \]
has the property, when multiplied by the matrix
\[ \begin{bmatrix}
    c_0\\
    \vdots\\
    c_{k+d}
\end{bmatrix},  \]
of yielding the homogeneous valuation of the polynomial $c(z)=\sum c_i z^i$ at the homogeneous point $(b,a)$, i.e., $\sum c_i b^i a^{k+d-i}$. Similarly, $\partial_a^l w(a,b)$ applied to the same matrix yields the homogeneous valuation at $(b,a)$ of the $l$-th derivative of $c(z)$. This immediately yields that, if $(b_j,a_j)$ is a root of $p(z)$ of multiplicity $l+1$, then the vectors
$w(a,b), \partial_a w(a,b), \dots, \partial_a^l w(a,b)$ belong to the null space of $\mathcal{T}_{k}(v^T)$ since the block columns of $\mathcal{T}_{k}(v^T)$ are the coefficients of the polynomials $v(z),zv(z),z^2v(z),\dots,z^k v(z)$, all of which have a root of multiplicity at least $l$ at $(b_j,a_j)$. This yields $d-m$ linearly independent vectors (linear independence follows from standard properties of, possibly confluent, Vandermonde matrices) in the null space of $\mathcal{T}_{k}(v^T)$. This shows that $\mathrm{rank} \,\mathcal{T}_{k}(v^T) \leq k+m+1$.
\end{proof}

Lemma \ref{lemma:rdrop_minindex} shows that the matrix $\mathcal{T}_k(v^T)$ experiences a rank drop when the degree of $v(x)$ drops or when $v(x)$ stops being a minimal basis. This suggests that the function
\begin{equation}
\label{eq:f_without_regulariz_for matrix_polynomials}
    f(v)= \norm*{\begin{bmatrix}
    A_0 & \ldots & A_k 
\end{bmatrix} \mathcal{T}_k(v^T)^T (\mathcal{T}_k(v^T)^T)^{\dagger}}_F
\end{equation}
can be discontinuous at the vectors $v$ for which either (or both) of these conditions holds. Thus, in this context, a regularization technique seems needed and the use of the penalty method is a possible way to overcome this issue.

In \cite{bora}, the authors propose a minimization technique for the function \eqref{eq:f_without_regulariz_for matrix_polynomials}, optimizing over $\mathbb{F}^{n}$, but they do not include a classical Tikhonov regularization strategy. Instead, they suggest an ad hoc strategy that consists in repeating the minimization process on a sequence of increasing degrees $d$ of the polynomial vector in the null space (corresponding to different submatrices of our $\mathcal{T}_k(v^T)$). An extensive comparison among our method and the algorithm proposed in \cite{bora} is provided in Subsection \ref{sec:betterthanbora}.

\section{Approximate GCD}\label{sec:approxgcd}

Let us consider two univariate polynomials with coefficients in $\C$,  say, $p(x)=\sum_{j=0}^m p_j x^j$ and $q(x)=\sum_{i=0}^n q_i x^i$, with $p_m \neq 0 \neq q_n$. A Greatest Common Divisor (GCD) of the pair $p(x),q(x)$ is a polynomial $g(x)=\sum_{k=0}^d g_k x^k$ such that (1) $g(x)$ divides both $p(x)$ and $q(x)$ (2) every polynomial satisfying property (1) divides $g(x)$. The GCD is unique up to multiplication by units, i.e., nonzero scalars.

The problem of computing the GCD between two polynomials is well-known to be ill-posed: Even if two polynomials $p$ and $q$ have nonconstant common factors, generically their perturbations have trivial GCD. Hence, it makes sense to consider the following ``approximate GCD" computational problem: Given two nonzero polynomials $p,q$ and a fixed degree $d \leq \min(\deg p, \deg q)$, find
\begin{equation} \label{mindistancegcd}
    f_d(p,q) = \min_{{\delta_p,\delta_q}} \, \norm{(\delta_p,\delta_q)} \subjectto \deg \operatorname{gcd}(p+\delta_p, q+\delta_q)\geq d.
\end{equation}
Here, $\delta_p$ and $\delta_q$ are polynomials of the same degree as $p$ and $q$ respectively, and 
\[
\norm{(\delta_p,\delta_q)} = \norm{\begin{bmatrix}
    (\delta_p)_0 & (\delta_p)_1 & \dots & (\delta_p)_m & (\delta_q)_0 & (\delta_q)_1 & \dots & (\delta_q)_n
\end{bmatrix}},
\]
i.e., the Euclidean norm of the vector of their coefficients. This problem has been studied extensively, together with variants where one imposes a maximum norm for $(\delta_p,\delta_q)$ and looks for the largest possible $d$; see, e.g.,\cite{Bini2010,KYZ,NAGASAKA} and the references therein for more information on the problem and its applications.

To see how GCD computations fit into our framework, let us start with a general algebraic result.
\begin{proposition} \label{prop:gcd}
    Let $R$ be a GCD domain with a multiplicative valuation $\nu$, i.e., a function $$z \in R \setminus \{0\} \mapsto \nu(z) \in \mathbb{N}=\{0,1,2,..\}$$ such that $\nu(0)=-\infty$ and $\nu(yz)=\nu(y)+\nu(z)$ $\forall$ $y,z \in R$. Let $p,q,g \in R$ with $g=\gcd(p,q)$. Then, $\nu(g)=\nu(p)-\nu(w)=\nu(q)-\nu(u)$ where $(u,w)\neq 0$ satisfy the equation $up+wq=0$ and have minimal valuations among  nonzero solutions.
\end{proposition}
\begin{proof}
$up+wq=0$ implies {that $\ell:= \mathrm{lcm}(p,q)$ divides $up$}. 
Thus, $\nu(q)-\nu(g) = \nu(\ell) -\nu(p) \leq \nu(u)$. Similarly, $\nu(w) \geq \nu(p)-\nu(g)$. To complete the proof, we exhibit the minimal solution $u=-\frac{q}{g}$, $w=\frac{p}{g}$.
\end{proof}
A multiplicative valuation as in Proposition \ref{prop:gcd} always exists when $R=\C[x]$, and it is usually called ``degree". We conclude that two polynomials $p(x),q(x)\in \C[x]$ have a GCD of degree $\geq d$ if and only if there are two nonzero polynomials $u(x)$ and $w(x)$ of degrees $\deg u(x) = \deg q(x) - d$, $\deg w(x) = \deg p(x) - d$ and such that $u(x)p(x) + w(x)q(x) = 0$. Essentially, $u(x)$ contains the roots of $q(x)$ that are not also roots of $p(x)$, and vice versa for $w(x)$.

If we use the notation of Section \ref{sec:Nearest_sing_matrix_poly} to turn polynomial multiplication into matrix-vector multiplication, we see that the problem of finding a solution $(u,w)$ to $up+wq=0$ is equivalent to seeking a nonzero vector
    \[
    \begin{bmatrix}
        \vvec u\\
        \vvec w
    \end{bmatrix} \in \ker 
     \begin{bmatrix}
        \mathcal{T}_{\deg u}(p) &
        \mathcal{T}_{\deg w}(q)
    \end{bmatrix}.
    \]
    We apply a block scaling so that norms are preserved: Define
    \[
    S_d(p,q) = \begin{bmatrix}
        \frac{1}{\sqrt{\deg(q)-d+1}} \mathcal{T}_{\deg(q)-d}(p) & 
        \frac{1}{\sqrt{\deg(p)-d+1}} \mathcal{T}_{\deg(p)-d}(q)
    \end{bmatrix}.
    \]
    Then, one sees that $\norm{(\delta_p,\delta_q)} = \norm{S_d(\delta_p,\delta_q)}_F$, and, in view of the previous discussion, $\deg \operatorname{gcd}(p+\delta_p, q+\delta_q) \geq d$ if and only if there is a nonzero vector in $\ker S_d(p+\delta_p,q+\delta_q)$.
    We can take $\mathcal{S}$ to be the (linear) space of matrices of the form $S_d(\cdot,\cdot)$, and then the problem~\eqref{mindistancegcd} can be recast into our framework of distance from singularity
    \[    
    f_d(p,q) = \min_{{\delta_p,\delta_q}} \, \norm{S_d(\delta_p,\delta_q)}_F \subjectto \text{$S_d(p,q) + S_d(\delta_p,\delta_q)$ is singular}.
    \]
    
It is worth remarking that we could speed up our computations (and potentially increase their accuracy) by exploiting the Toeplitz structure in the matrix $S_d$. In our numerical experiments, we do not currently use ad-hoc strategies that take this observation into account.

\section{Structured distance to instability} \label{sec:instability}

Given an open set $\Omega \subseteq \mathbb{C}$, we say that a square matrix is \emph{$\Omega$-stable} if all its eigenvalues belong to $\Omega$. Common examples are Hurwitz stability, with $\Omega = \{\lambda \in \mathbb{C} \colon \operatorname{Re} \lambda < 0\}$ the left half-plane, and Schur stability, with $\Omega = \{\lambda \in \mathbb{C} \colon \abs{\lambda} < 1\}$ the unit disc.

Another common problem in applications is the following~\cite{Bye88,HeWatson,Kre06,VL85}: Given an $\Omega$-stable matrix $A \in \mathbb{C}^{n\times n}$, compute its \emph{distance to $\Omega$-instability}, i.e.,
the norm of the smallest perturbation $\Delta$ such that $A+\Delta$ is not stable:
\[
\operatorname{dist}(A, \Omega^c) = \min_{{\Delta \in \mathbb{C}^{n\times n}}} \, \norm{\Delta}_F \subjectto \text{$A+\Delta$ has an eigenvalue $\lambda \in \Omega^c$},
\]
where we denote by $\Omega^c$ the complement set $\mathbb{C} \setminus \Omega$. Constraining the structure of the perturbation $\Delta$ is a challenging addition to this problem, studied for instance in~\cite{Sicilia}: It reflects the fact that real-life matrices from control systems often come with prescribed sparsity patterns or equal entries.

We recall a few facts on distances and projections; see~\cite[Section~2]{NP21} for more details. Given an open set $\Omega \subseteq \mathbb{C}$, we define the projection $\operatorname{proj}(z, \Omega^c)$ as the closest point (computed according to the Euclidean metric on the complex plane) to $z$ outside $\Omega$. This projection is unique for every $z$ except for a measure-zero set of singular points called the \emph{medial axis} of $\Omega^c$. Outside of the medial axis, the squared-distance function
\[
\operatorname{dist}^2(z, \Omega^c) = \abs{z - \operatorname{proj}(z, \Omega^c)}^2
\]
is a $\mathcal{C}^1$ function in $\mathbb{C} \equiv \mathbb{R}^2$ and has gradient (according to the Euclidean metric on $\mathbb{C} \equiv \mathbb{R}^2$)
\[
\operatorname{grad}_z \operatorname{dist}^2(z, \Omega^c) = 2(z - \operatorname{proj}(z, \Omega^c)).
\]

Following the approach in the previous sections, we can reformulate the distance to instability problem as
\[
\operatorname{dist}^2(A, \Omega^c) = \min_{{\Delta \in \mathcal{S}}} \, \norm{\Delta}_F^2 \subjectto (A+\Delta-\lambda I)v = 0 \text{ for some $v\in\mathcal{M}$ and $\lambda \in \Omega^c$},
\]
which is equivalent to the nested minimization problems
\begin{align*}
\operatorname{dist}^2(A, \Omega^c) &=  \min_{{ v\in\mathcal{M}}} \, f(v) \\
f(v) &= \min_{{\lambda\in\Omega^c}} \, f(v,\lambda) \\ 
f(v,\lambda) &= \min_{{\Delta \in \mathcal{S}}} \, \norm{\Delta}_F^2 \subjectto 
(A+\Delta-\lambda I)v = 0.
\end{align*}
We switch immediately to the relaxed formulation
\begin{subequations}
\begin{align}
\operatorname{dist}^2_\varepsilon(A, \Omega^c) &= \min_{{ v\in\mathcal{M}}} f_\varepsilon(v) \label{neareststab1} \\
f_\varepsilon(v) &= \min_{{ \lambda\in\Omega^c}} f_\varepsilon(v,\lambda) \label{neareststab2}\\
f_\varepsilon(v,\lambda) &= \min_{{ \Delta \in \mathcal{S}}} \norm{\Delta}^2_F + \varepsilon^{-1}\norm{(A + \Delta - \lambda I )v}^2.\label{neareststab3}
\end{align}
\end{subequations}

The problem~\eqref{neareststab3} is precisely the problem that we have solved in the last sections, except that $A$ is replaced by $A-\lambda I$. Hence,
\begin{equation} \label{fepsilonlambda}
f_\varepsilon(v,\lambda) = v^*(\lambda I-A)^* (MM^*+\varepsilon I)^{-1}(\lambda I - A)v.    
\end{equation}
We show that the minimizer $\lambda_*$ of~\eqref{neareststab2} can also be found explicitly, since $f_\varepsilon(v, \lambda)$ is a quadratic function of $\lambda$. Let us define the complex scalars
\begin{align*}
    a &= v^*(MM^*+\varepsilon I)^{-1}v \geq 0, &
    b &= v^*(MM^*+\varepsilon I)^{-1}Av \in \mathbb{C},\\
    c &= v^*A^*(MM^*+\varepsilon I)^{-1}Av  \geq 0, &
    \lambda_0 &= \frac{b}{a}.
\end{align*}
\begin{theorem} \label{thm:minimizerproject}
    {Let $f_\varepsilon(v,\lambda)$ be defined in~\eqref{fepsilonlambda}. Then,} the minimizer $\lambda_* = \arg \min_{\lambda\in\Omega^c} f_\varepsilon(v,\lambda)$ is given by 
    $\lambda_* = \operatorname{proj}_{\Omega^c} (\lambda_0)$.
\end{theorem}
\begin{proof}
We have
\begin{align*}
    f_\varepsilon(v,\lambda) &= (\overline{\lambda}v^* - A^*)(MM^*+\varepsilon I)^{-1}(v\lambda - Av)\\
    &= \overline{\lambda} a \lambda - \overline{\lambda}b - \overline{b}\lambda + c\\
    &= (\overline{\lambda - a^{-1}b})a(\lambda - a^{-1}b) - \overline{b} a^{-1} b+c\\
    &= a \abs{\lambda - \lambda_0}^2 + c - a^{-1}\abs{b}^2.
\end{align*}
The result follows from noting that the level curves of $f_\varepsilon(v, \lambda)$, plotted in the complex plane as a function of $\lambda$, are circles centered in $\lambda_0$.
\end{proof}
 Hence, to evaluate $f_\varepsilon(v)$ for this problem, we can use Algorithm~\ref{algo:f_unstable}.
\begin{algorithm}
\caption{Computing $f_\varepsilon(v)$ for the nearest unstable matrix problem} \label{algo:f_unstable}
    $M \gets [P^{(1)}v, \dots, P^{(p)}v]$ \;
    $\lambda_0 \gets \frac{v^*(MM^*+\varepsilon I)^{-1}Av}{v^*(MM^*+\varepsilon I)^{-1}v}$\;
    $\lambda_* \gets \operatorname{proj}_{\Omega^c}(\lambda_0)$\;
    $f_\varepsilon(v) \gets v^*(\lambda_* I-A)^* (MM^*+\varepsilon I)^{-1}(\lambda_* I - A)v$\;
\end{algorithm}
Computing the gradient of this function looks like a difficult task at first, but it is simplified by the following observation: If $\operatorname{proj}(\lambda_0, \mathbb{C}) = \zeta$, for a certain $\zeta \in \Omega^c$, then the functions 
$\operatorname{dist}^2(\lambda_0, \Omega^c)$ and $\abs{\lambda_0 - \zeta}^2$ (where $\zeta$ is considered as a constant) have the same first-order expansion in $\lambda_0$, since they have the same value and the same gradient $2(\lambda_0 - \zeta)$.

In particular, the gradient of
\[
f_\varepsilon(v) = f_\varepsilon(v, \lambda_*) = a \operatorname{dist}^2(\lambda_0, \Omega^c) + c - a^{-1}\abs{b}^2 = a \abs{\lambda_0 - \lambda_*}^2 + c - a^{-1}\abs{b}^2
\]
with respect to $v$ is the same as the gradient of
\[
f_\varepsilon(v, \lambda_*) = a \abs{\lambda_0 - \zeta}^2 + c - a^{-1}\abs{b}^2
\]
with the constant $\zeta = \lambda_*$, i.e., we can regard $\lambda_*$ as a constant in our gradient computations and ignore the contribution from $\frac{d}{d t} \lambda_*$. Hence, we can repeat the computation of the previous section with $A-\lambda_* I$ in place of $A$, to get
\[
\nabla f_\varepsilon(v) = -2(A-\lambda_* I +\Delta_*)^*(MM^*+\varepsilon I)^{-1}(\lambda_* I - A)v.
\]
\begin{remark}
    The argument above does not extend to the computation of the Hessian, unfortunately, since the Hessians of $\operatorname{dist}^2(\lambda_0, \Omega^c)$ and $\abs{\lambda_0 - \zeta}^2$ are different. This fact suggests using first-order algorithms for this problem, since they can be implemented relying only on the function  $z\mapsto \operatorname{proj}(z, \Omega^c)$, without the need for additional information such as its Hessian.
\end{remark}
\begin{remark}
    For this problem, the field $\mathbb{F}$ must be $\mathbb{C}$: Even when $A$ has real entries, $\lambda_*$ and $v_*$ may be complex; hence, we must look for the minimum in the complex field. One can also consider the problem of finding the nearest unstable $A+\Delta$, with $\mathcal{S}\subseteq \mathbb{R}^{n\times n}$ equal to the set of \emph{real} linear combinations of a basis set $P^{(1)}, \dots, P^{(p)}$. We can proceed similarly, imposing $\delta\in\mathbb{R}^{p}$ and converting the complex linear system $M\delta = r$ into the real one
    \[
        \min_{{\delta \in \mathbb{R}^p}}\, \norm{\delta}^2 {{}\subjectto{}}
        \begin{bmatrix}
            \operatorname{Re} M\\
            \operatorname{Im} M\\
        \end{bmatrix}\delta = 
        \begin{bmatrix}
            \operatorname{Re} r\\
            \operatorname{Im} r
        \end{bmatrix}.
    \]
    However, an analogue of Theorem~\ref{thm:minimizerproject} does not hold: In this case, finding the optimal $\lambda_*$ is a general convex quadratic minimization problem over $\mathbb{R}^2$, and its solution is not given by the projection in the Euclidean distance.
\end{remark}

\section{Nearest matrix with prescribed nullity}
\label{sec:nearest_matrix_pescribed_nullity}
In this section, we study another nearness problem that we can solve with a similar strategy. Given $A\in\mathbb{F}^{m\times n}$, $m\geq n$, we wish to find the nearest matrix $A+\Delta$ with $\rank (A+\Delta) \leq n-\ell$, i.e., nullity at least $\ell$. When $\ell=1$, this problem reduces to the one of the nearest singular matrix that we have already treated. 

A matrix $A+\Delta$ has nullity at least $\ell$ if and only if there exists a $V\in\mathbb{F}^{n\times \ell}$, with orthonormal columns, such that $(A+\Delta)V = 0$. Hence, we define
\[
f(V) = \min_{{\Delta\in\mathcal{S} }} \, \norm{\Delta}^2 \subjectto (A+\Delta)V = 0.
\]
This function is invariant under changes of bases $V \mapsto VQ$, with $Q$ unitary, so we can optimize over the Grassmann manifold \cite{EdelmanGrassmann} $G_{n,\ell}(\mathbb{F})$. Here, $V$ is an orthonormal representative of an element of $G_{n,\ell}(\mathbb{F})$.

Set
\[
V = \begin{bmatrix}
    v_1 & v_2 & \dots & v_\ell
\end{bmatrix};
\]
then $(A+\Delta)V = 0$ is equivalent to $(A+\Delta)v_i = M(v_i)\delta - r(v_i)=0$, for all $i=1,2,\dots,\ell$. Hence, we can rewrite the problem as a minimization problem in $\delta$ as
\[
f(V) = \min_{{ \delta \in \mathbb{F}^{p} }} \, \norm{\delta}^2 \subjectto M(v_i)\delta = r(v_i),\, i=1,2,\dots,\ell,
\]
or equivalently
\[
f(V) = \min_{{\delta\in\mathbb{F}^p}}\, \norm{\delta}^2 \subjectto M(V)\delta = r(V), \quad
M(V) = 
\begin{bmatrix}
    M(v_1)\\
    M(v_2)\\
    \vdots\\
    M(v_\ell)
\end{bmatrix}, \quad r(V) = \begin{bmatrix}
    r(v_1)\\
    r(v_2)\\
    \vdots\\
    r(v_\ell)
\end{bmatrix}.
\]
The formulae to compute $f$ in the previous sections apply almost identically to this problem, with the only caveat that the sizes might be different: While in the previous cases we expected that $M$ has more rows than columns, in this case the opposite is more typical. This is relevant in the formulae for the pseudoinverse, which becomes $M^\dagger = (M^*M)^{-1}M^*$ in the full-column case; however in the regularized version we still have $M^*(M M^*+\varepsilon I)^{-1} = (M^* M+\varepsilon I)^{-1}M^*$.

\section{Numerical experiments}\label{sec:numericalexperiments}

In this section, we provide details on the numerical implementation of our method, and present extensive numerical experiments. In Subsection \ref{subsec:adaptive strategy} we describe our heuristic strategy for the choice of the regularization parameter $\varepsilon$ in the penalty method, stated in Subsection \ref{subsec:penalty-method-augmented-Lagrangian}. The remaining part of this section is devoted to numerical experiments for different applications of our method. We test our method using as benchmark pre-existing algorithms. The Matlab code for the implementation is freely available at \url{https://github.com/fph/RiemannOracle}. All running times were measured using Matlab R2023a and Manopt 7.1 on an Intel Core i5-12600K.

\subsection{Adaptive choice for the regularization parameter}
\label{subsec:adaptive strategy}

Depending on the method used for minimization, we may use a different strategy to update the parameter $\varepsilon$ at each iteration. In particular, in the penalty method, as the parameter $\varepsilon$ must tend to $0$, we need to generate a (finite) decreasing sequence $(\varepsilon_k)_{k \leq K}$ and solve minimization problems for smaller and smaller values of $\varepsilon_k$, $k=0,1,\dots,K$. It is common practice \cite{NoceWrig} to use the solution at the step $k$, that is,
\[
v_k^* = \arg \min_{v \in \mathcal{M}} f_{\varepsilon_k}(v),
\]
as a starting point for minimizing $f_{\varepsilon_{k+1}}(v)$. Clearly, it is desirable that $K$ is as small as possible, for reasons of efficiency. Nevertheless, decreasing the regularization parameter too fast may lead to spurious results. We provide a heuristic adaptive strategy for choosing the value $\varepsilon_{k+1}$ from $\varepsilon_k$, described in Algorithm \ref{algo:adaptive} and motivated by the following observations. 

\begin{figure}[h!]
    \centering
    \includegraphics[width=0.92\textwidth]{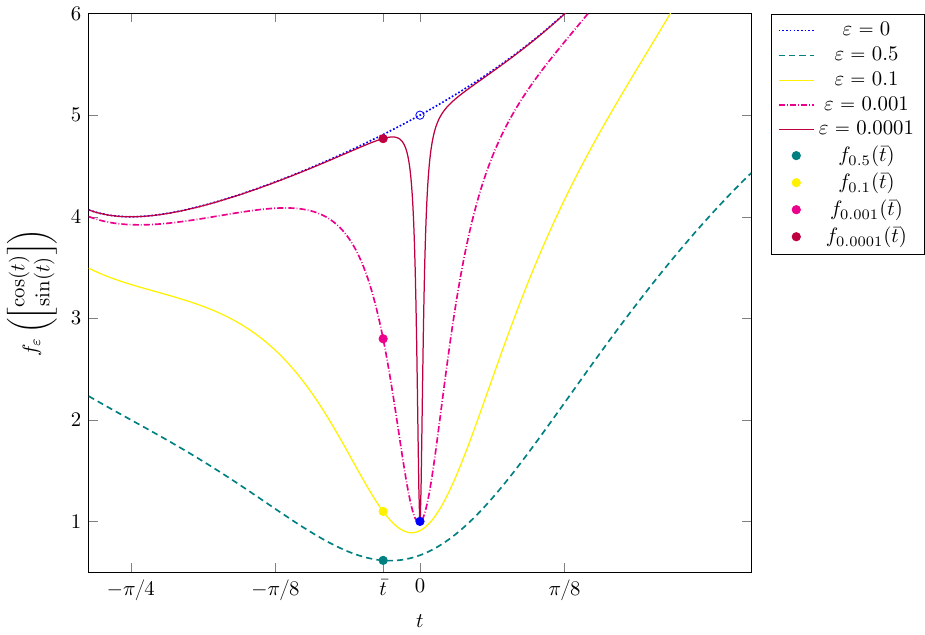}
    \caption{A zoom-in of Figure \ref{fig1}.}
    \label{fig:Zoom on penalty}
\end{figure}

In Figure \ref{fig:Zoom on penalty}, we provide a zoom-in of Figure \ref{fig1}. Having fixed a point $\bar{t}$, we consider the points $f_{\varepsilon}(\bar{t})$, where $\varepsilon= 0.5$ and $\varepsilon=10^{-j}$, for $j=1,2,3$. In this example, the geometric decrease of $\varepsilon$ induces a relatively large difference between the values of $f_{\varepsilon}(\bar{t})$. This may lead to possible drawbacks. Firstly, the point $v_{k}^*$ may be not a suitable starting point for iteration $k+1$. Secondly, it is possible that the minimization procedure converges to a different local minimum. For instance, in the example provided in Figure \ref{fig:Zoom on penalty}, this happens if, after solving the subproblem for $\varepsilon_k=0.5$, we consider $\varepsilon_{k+1} = 0.0001$ and use $\bar{t}$ as starting point for minimizing $f_{0.001}(v)$. 

Heuristically, we observed that these issues can be avoided or greatly reduced {by employing a backtracking strategy: we aim to set $\varepsilon_{k+1} =  \varepsilon_{k}\mu$, for a certain reduction factor $\mu<1$. We try initially with a very small value $\mu$; however, if this choice of $\mu$ leads to a value of $f_{\varepsilon_{k+1}}(v_k^*)$ that is much larger than $f_{\varepsilon_{k}}(v_k^*)$, then we ``backtrack'' and try a slightly larger $\mu$.}

We used this technique when solving the problem stated in Subsection \ref{subsec:unstructured_distance_polynomial}. We report these experiments in Subsection \ref{sec:betterthanbora}. {In our implementation in \url{https://github.com/fph/RiemannOracle}, we chose $[\mu_{\min},\mu_{\max}] = [0.01,0.95]$, $\rho=1.1$, $C_f=2.5$.}

\begin{algorithm}[h!]
\caption{Adaptive choice of $\varepsilon_{{k+1}}$} \label{algo:adaptive}
\SetKwInOut{Input}{Input}
\SetKwInOut{Output}{Output}
\Input{Starting value $\varepsilon_k$, {range $[\mu_{\min},\mu_{\max}]$ with $\mu_{\min}<\mu_{\max}<1$,} {backtracking factor $\rho>1$, constant $C_f>1$}}
\Output{New choice $\varepsilon_{k+1}$}
$v_k \gets \arg \min_{v\in\mathcal{M}} f_{\varepsilon_k}(v)$ \\
$f_k \gets f_{\varepsilon_{k}}(v_k)$\\
{$\mu \gets \mu_{\min}$}\\
\While{${f_{\varepsilon_{k}\mu}(v_k)} > C_f f_k$}{
     $\mathstrut{\mu \gets \mu\rho } $ \\
    \If{${\mu> \mu_{\max} } $}{
     \strut{\textbf{break}}\\
    }
}
{$\varepsilon_{k+1} \gets \varepsilon_k \mu$} \\
\end{algorithm}

\subsection{Nearest sparse singular matrix: A comparison with \cite{Sicilia}}\label{sec:betterthansicilia}
\begin{figure}
    \centering
    \begin{minipage}{0.49\linewidth}
        \includegraphics[width=\linewidth]{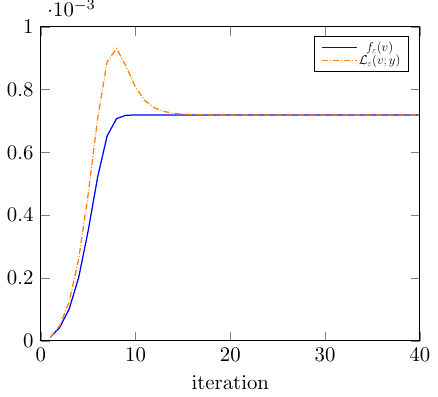}
    \end{minipage}\hfill
    \begin{minipage}{0.49\linewidth}
        \includegraphics[width=\linewidth]{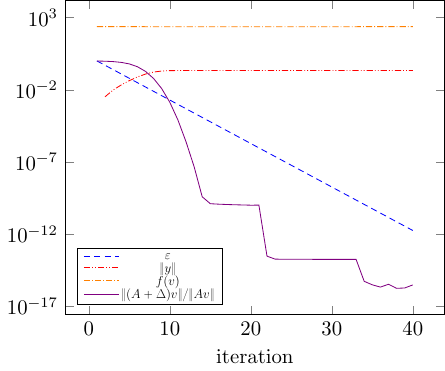}
    \end{minipage}\\[1em]
    \begin{minipage}{0.49\linewidth}
        \includegraphics[width=\linewidth]{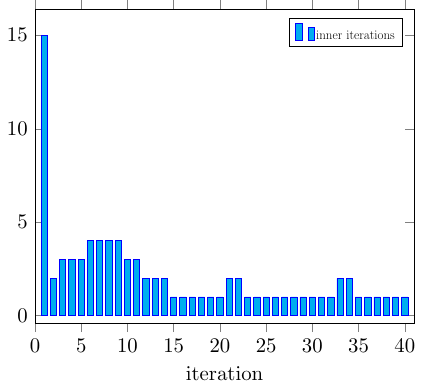}
    \end{minipage}\hfill
    \begin{minipage}{0.49\linewidth}
        \includegraphics[width=\linewidth]{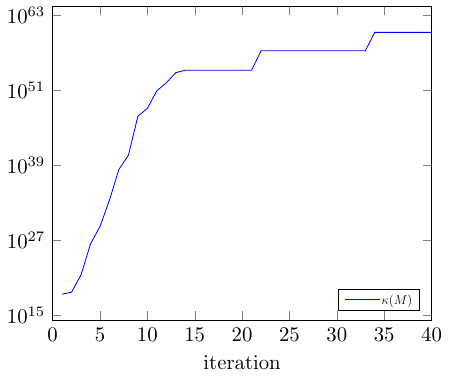}
    \end{minipage}
    \caption{Convergence plots for the nearest sparse singular matrix to \texttt{orani678}.}
    \label{fig:orani}
\end{figure}

In this experiment, given a sparse matrix $A$, we aim to compute the nearest singular matrix with the same zero pattern as $A$. We use the optimized formulae in Section~\ref{sec:sparse}. We take $A$ to be the sparse matrix \verb!orani678! from the Suitesparse Matrix Collection ({whose dimension is $n=2529$}, {with} $90158$ nonzero elements). We use the augmented Lagrangian method, decreasing $\varepsilon$ by a fixed factor $2$ at each iteration. The results are shown in Figure~\ref{fig:orani}.

The matrix $M$ is highly ill-conditioned at every step of the algorithm, making this problem a particularly challenging benchmark. Nevertheless, the value of $f_\varepsilon(v)$ converges quickly to $7.1894\times 10^{-4}$, corresponding to a structured distance from singularity of $2.6813\times 10^{-2}$.

It is tricky to evaluate the non-regularized function $f(v_*)$ at the minimizer vector $v_*$: Indeed, the second subplot in Figure~\ref{fig:orani} shows that the computed value of $f(v)$ stays at a value of about $240$, unlike $f_\varepsilon(v)$. This phenomenon happens because in the last iterations $M$ has $240$ singular values (out of $2529$) smaller than the machine precision, some as small as $10^{-57}$: Using the notation in Section~\ref{sec:sparse}, both $r_i$ and $d_i = \sigma_i^2$ are zero only up to numerical errors of the order of the machine precision, for large $i$, and the computation of $f$ is sensitive to their exact values.

At the exact minimizer, we expect that $d_i = r_i = 0$. If we enforce this property, we can produce a nearby vector $v_\mathrm{reg}$ for which the minimum $f(v_\mathrm{reg}) = 7.1894\times 10^{-4}$ is achieved even for $\varepsilon=0$: We set
\[
(r_\mathrm{reg})_i = \begin{cases}
    r_i & d_i > 10^{-16},\\
    0 & \text{otherwise},
\end{cases}
\]
and compute $v_\mathrm{reg} = -A^{-1}r_\mathrm{reg}$. This correction ensures that $r\in \operatorname{Im} M$ exactly, even in the limit when some entries $d_i$ become zero.

It is interesting to compare the values computed by Riemann-Oracle with those reported in~\cite[Section~7.1]{Sicilia}, which deals with similar problems. The structure of their computation is different: Their main algorithm solves the dual problem
\[
\phi(\varepsilon) = \min_{{\Delta\in\mathcal{S}}} \abs{\lambda}^2 \subjectto \text{$A+\Delta$ has eigenvalue $\lambda$,} \, \norm{\Delta}_F = \varepsilon,
\]
and then an outer minimization routine on $\varepsilon$ must be used to find the smallest $\varepsilon$ such that $\phi(\varepsilon)=0$.

We can obtain directly comparable results by solving a distance-to-instability problem, in the framework of Section~\ref{sec:instability}, setting $\Omega = \{z\in\mathbb{C} \colon \abs{z}\geq \phi(\varepsilon)^{1/2} \}$. In this way, we can compute the inverse of the function $\phi$: Given a target magnitude $\phi(\varepsilon)^{1/2}$ for the smallest eigenvalue of $A+\Delta$, we obtain the smallest perturbation norm $\min_v f(v) = \norm{\Delta}_F^2$ that achieves it. We report our results in Table \ref{table:compareSicilia}. The
agreement is not perfect up to the last reported digit, but we note the disclaimer in~\cite[Section~7.1]{Sicilia} stating that $\phi$ is computed inexactly; we can view this experiment as a confirmation that both methods produce valid results. Observe also the duality between the two frameworks: In~\cite{Sicilia}, intermediate values of $\varepsilon$ and $\phi(\varepsilon)$ are computed naturally as intermediate results to solve the distance to singularity problem; on the other hand, it is relatively easy for Riemann-Oracle to compute $\phi^{-1}(0)$, i.e., the distance to singularity, and a more sophisticated (and slower) variant is required to compute $\phi^{-1}(\varepsilon)$ for $\varepsilon > 0$.

The authors of \cite{Sicilia} report a CPU time of more than 500 seconds per outer iteration, with at least 10 outer iterations. While comparing running times on different machines is of course very imprecise, it is worth mentioning that Riemann-Oracle provides an output in less than a minute; this appears to be a speedup of several orders of magnitude.
\begin{table}[h!]
\centering
\begin{tabular}{cc}
    \toprule
     $\varepsilon = \norm{\Delta}_F$ & $\phi(\varepsilon) = \abs{\lambda_{\min}(A+\Delta)}^2$ \\
     \midrule
     $\underline{0.0000000}$ & $1.1019564 \cdot 10^{-2}$ \\
     $0.0176409$ & $9.5284061 \cdot 10^{-4}$ \\
     $0.0219541$ & $2.5263758 \cdot 10^{-4}$\\
     $0.0243116$ & $6.5050153 \cdot 10^{-5}$ \\
     $0.02\underline{64097}$ & $1.6503282 \cdot 10^{-6}$ \\
     $0.0261739$ & $4.1561289  \cdot 10^{-6}$ \\
     $0.0264923$ & $1.0428313  \cdot 10^{-6}$ \\
     $0.0266524$ & $2.6118300  \cdot 10^{-7}$ \\
     $0.0267327$ & $6.5355110  \cdot 10^{-8}$  \\
     $0.0267\underline{822}$ & $9.6346192  \cdot 10^{-9}$ \\
     $0.026\underline{8089}$ & $1.7293467  \cdot 10^{-10}$ \\
     \bottomrule
\end{tabular} 
\caption{Minimizers obtained by running the distance-to-instability algorithm in Section~\ref{sec:instability} with $\Omega =  \{z\in\mathbb{C} \colon \abs{z}\geq r \}$ and various values of $r$. The right column reports the least squared norm of the eigenvalues of the computed matrix. Underlined digits differ from \cite[Table 7.2]{Sicilia}.  We believe that our computations are likely more reliable for the first row (since $A$ already has an eigenvalue with a squared norm $1.1019564\cdot 10^{-2}$) and the sixth row (since in \cite{Sicilia} the function behaves nonmonotonically). For the two bottom rows, our values are slightly larger and we cannot say for sure that they are more accurate, but the authors of \cite{Sicilia} make a comment that their $\phi$ is computed with low accuracy.}\label{table:compareSicilia}
\end{table}

The software~\cite{MarU14software}, which implements a variant of the algorithm in~\cite{UseM14} to find the nearest singular matrix with a prescribed structure, can also be used with the sparsity structure. However, that code is not designed to deal with large, sparse matrices, and in particular it does not use the specialized formulae for the objective function described in Section~\ref{sec:sparse}. As a result, our attempts to use it with the \texttt{orani678} matrix failed with out-of-memory errors, probably due to the attempt to allocate $M$ as a full matrix.

Even if it could tackle such a large input, we conjecture that the approach in the \texttt{slra} code is likely to fail on this example, because  the global minimum is achieved at a discontinuity, hence involving the same issues discussed in Section~\ref{sec:need_regularization} and Remark~\ref{rem:smallcounterex}.

\subsection{Nearest singular matrix polynomial: A comparison with \cite{bora}} \label{sec:betterthanbora}
    
In this section, we focus on the nearest singular matrix polynomial to a given one of degree $\geq 2$; for the degree $1$ case, to our knowledge, the currently best available algorithm is \cite{DNN24} which already uses similar techniques. We compare the method of this paper with the algorithm presented in \cite{bora}\footnote{For the comparison, we have used Matlab codes kindly provided by the authors of \cite{bora}.}, which at the time of writing this paper is the state of the art for this problem for degree $2$ or higher. We do this by performing statistical experiments with randomly generated matrix polynomials. For a fixed size
$n$, we generate $n \times n$ matrix polynomials $A_0 + \lambda A_1 + \ldots +  \lambda^d A_d$ such that the real and imaginary part of each element in $A_i$, for $i=0,\ldots,d$ is drawn from the normal distribution $\mathcal{N}(0,1)$. 

The method proposed in \cite{bora} attempts to solve the unregularized problem \eqref{eq:unregular_polynomial_problem} by using a pseudoinverse in a similar manner as described in Theorem \ref{thm:fv}. The main difference between our method and the one in \cite{bora} is our use of regularization, as explained in Section \ref{sec:need_regularization}. Another significant difference is that we restrict to optimizing over the unit sphere, whereas the method in \cite{bora} does not normalize the input and optimizes over the Euclidean space of one dimension higher. 

The residual
\begin{align*} 
\eta(v, \Delta):=\norm*{\begin{bmatrix}
    \Delta_0+A_0 & \ldots & \Delta_k + A_k
\end{bmatrix}  \mathcal{T}_k(v^T)^T}_2
\end{align*}
measures the extent to which the constraint of the optimization problem \eqref{eq:unregular_polynomial_problem} is satisfied. Note that, in the case of the method presented in this paper, this residual approaches zero as the regularization parameter approaches zero. In the case of the method of \cite{bora}, on the other hand, this residual depends on the value of the tolerance used in the practical implementation of the pseudoinverse. In order to provide a fair comparison of the running times, we impose that the resulting residual of our method is less than or equal to that of the method in \cite{bora}. We do this as follows: We first run the algorithm of \cite{bora}, and measure the residual $\eta(\hat v, \hat \Delta)$, where $\hat v$ and $\hat \Delta$ are the minimizer and the corresponding perturbation, respectively, obtained by the algorithm in \cite{bora}. We stop the Riemann-Oracle as soon as $\eta(\tilde v, \tilde \Delta) \leq \eta (\hat v, \hat \Delta)$, where $\tilde v, \tilde \Delta$ are the minimizer and the corresponding perturbation, respectively, obtained by our approach.\footnote{In addition to the value of the stopping criterion, various other parameters need to be specified when running Manopt. The file \texttt{example\_polynomial.m} in the \href{https://github.com/fph/RiemannOracle}{GitHub} repository gives the complete list of the choices of the parameters that were used in {the experiments of this subsection}. These are stored in the \texttt{options} structure.} 

\begin{figure}[h!]
    \centering
    \includegraphics[width=\linewidth]{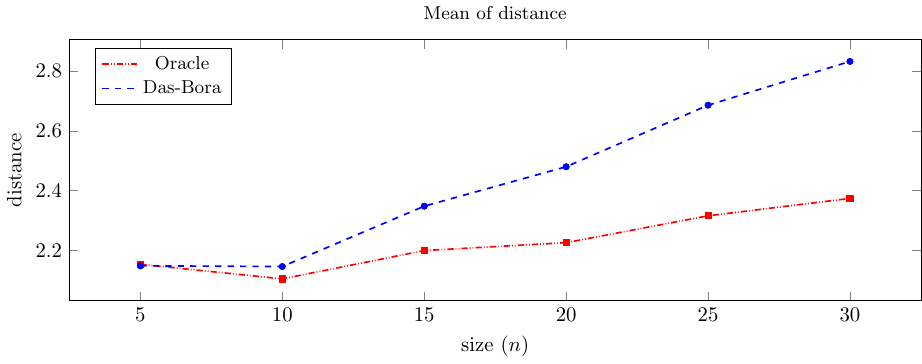}\\[1em]
    \includegraphics[width=\linewidth]{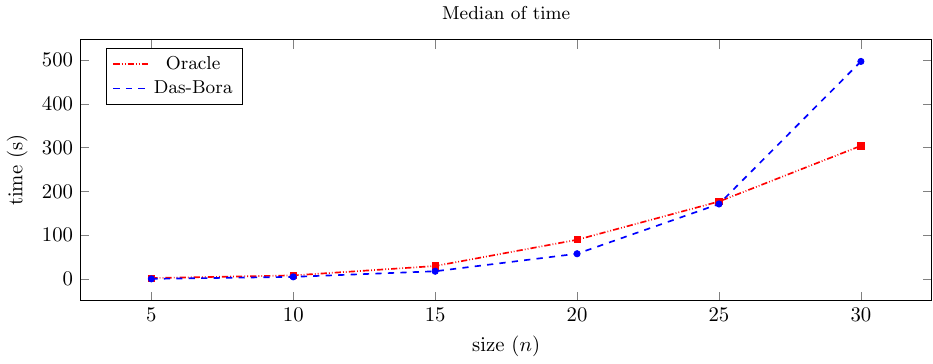}
    \caption{Comparison between Riemann-Oracle (denoted Oracle in the picture) and the Das-Bora algorithm \cite{bora} for $d=2$ and $n \in \{5,10,15,20,25,30\}$.}
    \label{fig:bora_comparison}
\end{figure}

\begin{figure}[h!]
    \centering
\includegraphics{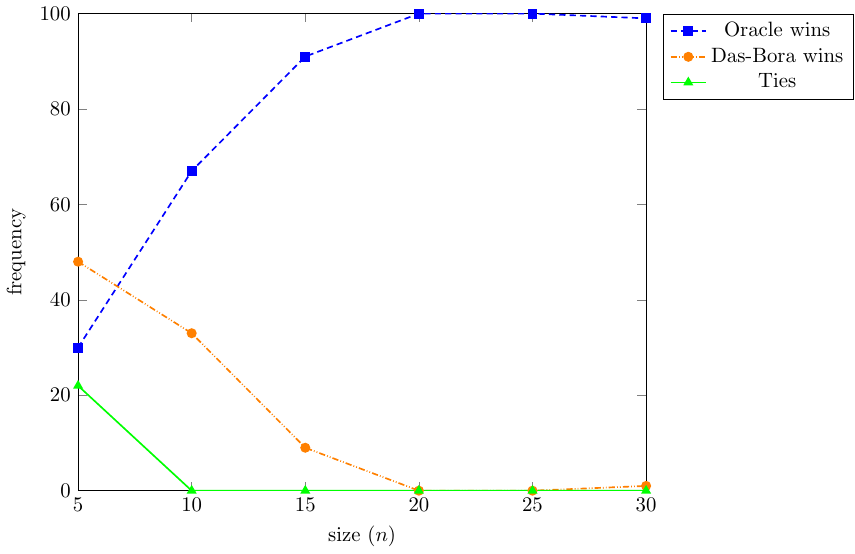}
    \caption{Comparison between the method of this paper (Oracle) and the Das-Bora algorithm \cite{bora} for $d=2$ and $n \in \{5,10,15,20,25,30\}$. The picture reports the relative frequency of which method yielded a better solution (or of ties).  A comparison is considered to be a tie if the solutions differ by at most $10^{-8}$.}
    \label{fig:bora_winloss}
\end{figure}

We design two different numerical experiments with the goal of comparing the method in \cite{bora} and Riemann-Oracle, implemented as detailed in Subsection \ref{subsec:unstructured_distance_polynomial}. Note that \cite{bora} cannot deal with additional structures, but, on the other hand, it can compute distances in the spectral norm as well. We compare the two algorithms on the intersection of their possibilities, that is, on computing the unstructured distance to singularity of a matrix polynomial in the Frobenius norm. In the first experiment, we consider randomly generated matrix polynomials of fixed degree $d=2$; the quadratic case $d=2$ is arguably the most relevant for practical applications \cite{TM}. For each $n \in \{5,10,15,20,25, 30\}$, we construct $100$ quadratic matrix polynomials. Figure \ref{fig:bora_comparison} shows a comparison between the two methods, both for the values of the approximate distances and the running times (in seconds). Figure \ref{fig:bora_winloss} shows how frequently each method yields a lower value of the objective function $\norm{\Delta}_F$. These figures show that as the size increases, the Riemann-Oracle approach outperforms its competitor, increasingly more impressively. In particular, for $n=30$, Riemann-Oracle finds a better solution than the method in \cite{bora} $99$ times out of $100$, while the average distance computed by our approach is approximately $16 \%$ smaller than the one found by \cite{bora}; this is a clear improvement to the state of the art.

\begin{figure}[h]
    \centering
    \includegraphics[width=\linewidth]{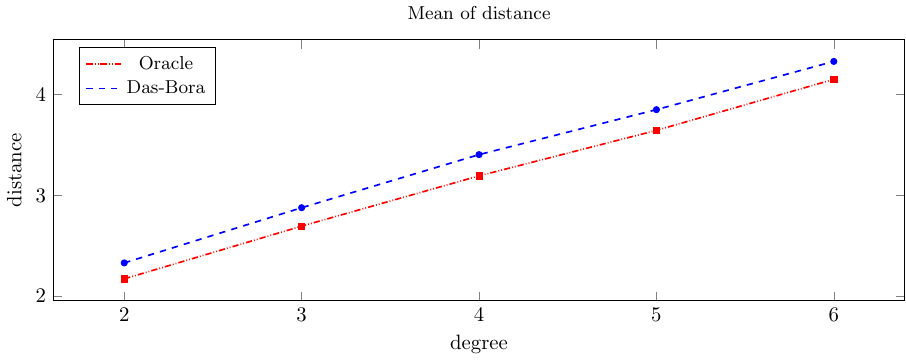}\\[1em]
    \includegraphics[width=\linewidth]{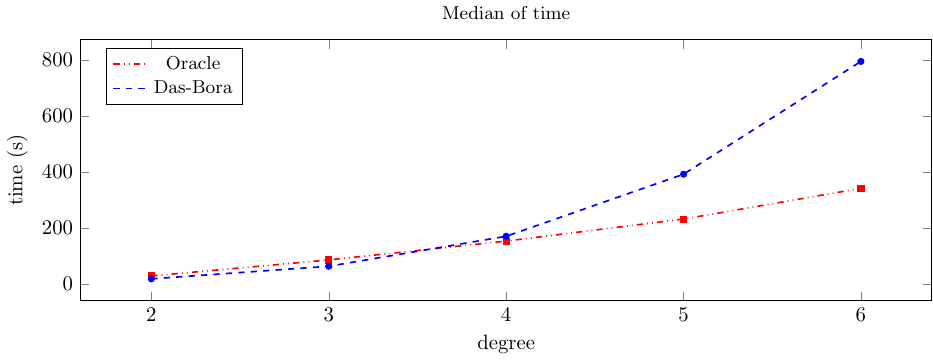}
\caption{Comparison between Riemann-Oracle (denoted Oracle in the picture) and the Das-Bora algorithm \cite{bora} for $d \in \left\lbrace 2, 3, 4 , 5 , 6 \right\rbrace$ and $n=15$. }\label{fig:bora_comparison_degree}
\end{figure}

\begin{figure}[h!]
    \centering
\includegraphics{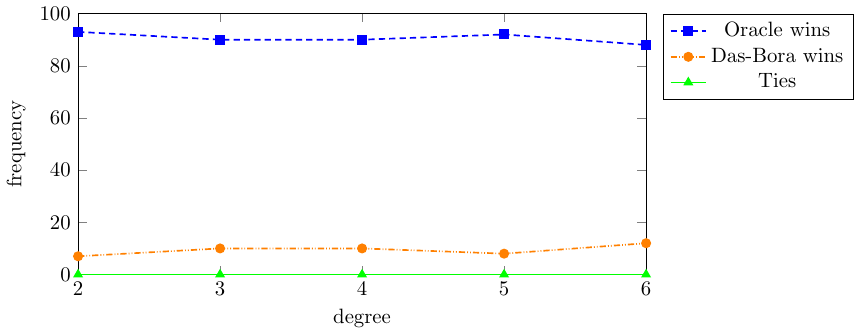}
\caption{Comparison between Riemann-Oracle and the Das-Bora algorithm \cite{bora} for $d \in \left\lbrace 2, 3, 4 , 5 , 6 \right\rbrace$ and $n=15$. The picture reports the relative frequency of which method yielded a better solution (or of ties). A comparison is considered to be a tie if the solutions differ by at most $10^{-8}$.}\label{fig:bora_winloss_degree}
\end{figure}

Next, we perform a similar experiment, but varying the degree of the matrix polynomial instead of its size. More specifically, we set $n=15$ and let the degree vary in the range $d \in \left\lbrace 2, 3, 4 , 5 , 6 \right\rbrace$. The results are shown in Figures \ref{fig:bora_comparison_degree} and \ref{fig:bora_winloss_degree}. For all degrees, the quality of the output given by Riemann-Oracle is slightly better than its competitor. In terms of running time, Riemann-Oracle scales significantly better for higher degrees. In particular, for $d=6$, the median running time of Riemann-Oracle is approximately one third of that of the Das-Bora algorithm.

\subsection{Approximate GCD: A comparison with \cite{NAGASAKA}}\label{sec:herewelose}

We apply the approach in Section~\ref{sec:approxgcd} to finding approximate GCDs of given degree $d$ of
\begin{equation} \label{boito861}
    p(x) = \beta(x^3 + 3x - 1)  (x-1)^k, \quad q(x) = \gamma\frac{d}{dx}p(x), 
\end{equation}
and of
\begin{equation} \label{boito811}
p(x) = \beta\prod_{j=1}^{10} (x-\alpha_j), \quad q(x) = \gamma\prod_{j=1}^{10} (x-\alpha_j + 10^{-j}), \quad \alpha_j = (-1)^j\tfrac{j}{2},
\end{equation}
with the normalization coefficients $\beta$ and $\gamma$ chosen so that the vectors of coefficients of $p$ and $q$ have Euclidean norm 1. These specific pairs of polynomials have been considered in various papers in the GCD literature~\cite{Bini2010,NAGASAKA,Zeng} as challenging test problems.

We apply the augmented Lagrangian method with a very small tolerance $\norm{\operatorname{grad} f_{\varepsilon,y}(v)} < 10^{-12}$, a maximum of 5000 steps in each inner iteration, and 80 outer iterations in which $\varepsilon$ is decreased by a constant factor $0.7$, starting from $1$. The results obtained by our algorithm on the two problems are reported in Table~\ref{tbl:gcd}. Note that we compute the polynomials $u,w,g$ explicitly and set $\delta_p = p - gu$, $\delta_q = q - gw$, to ensure that the computed distances truly correspond to a GCD of a given degree. The values obtained by Riemann-Oracle match the best ones reported in~\cite[Table~4]{NAGASAKA}, apart from the last two rows in the rightmost table: $d=5$, for which a slightly smaller value $4.487\cdot 10^{-9}$ is obtained by some of the specialized algorithms; and $d=4$, for which our algorithm converges to a local minimum corresponding to the cofactors of the approximate GCD with $d=6$, multiplied by two spurious degree-1 factors. Hence, a GCD polynomial $g$ of degree 4 cannot be extracted from those cofactors.

\begin{table}[]
\centering
    \begin{tabular}[t]{cc}
    \begin{tabular}[t]{ccc}
    \toprule
    \multicolumn{3}{c}{Approximate GCDs of~\eqref{boito861}}\\
    $k$ & $d$ & $\norm{(p-gu, q-gw)}$\\
    \midrule
    15 & 15 & $7.6160\cdot 10^{-5}$ \\
    25 & 25 & $7.8733\cdot 10^{-6}$ \\
    35 & 36 & $6.1775\cdot 10^{-5}$ \\
    45 & 46 & $2.9846\cdot 10^{-5}$ \\
    \bottomrule
    \end{tabular}
    &
        \begin{tabular}[t]{cc}
    \toprule
    \multicolumn{2}{c}{Approximate GCDs of~\eqref{boito811}}\\
    $d$ & $\norm{(p-gu, q-gw)}$\\
    \midrule
        9 & $3.9964\cdot 10^{-3}$\\
        8 & $1.7288\cdot 10^{-4}$\\
        7 & $7.0890\cdot 10^{-6}$\\
        6 & $1.8293\cdot 10^{-7}$\\
        5 & $4.4\underline{913}\cdot 10^{-9}$\\
        4 & \underline{$1.8293\cdot 10^{-7}$}\\
    \bottomrule
    \end{tabular}
    \end{tabular}
    \caption{Results obtained for the approximate GCD problem~\eqref{mindistancegcd} on the polynomials~\eqref{boito861} and~\eqref{boito811}, with varying degrees $k$ and $d$. The {underlined digits are those that differ from the results obtained by the best method in~\cite{NAGASAKA}. Further discussion is in the main text.}}
    \label{tbl:gcd}
\end{table}

We remark that we applied Riemann-Oracle without any post-processing of the solutions, unlike in~\cite{NAGASAKA} in which post-processing is considered as an integral part of the algorithm. The results confirm that our algorithm matches the accuracy of the best specialized algorithms, even without post-processing, apart from the more challenging cases in which the objective function takes a value smaller than the machine precision (note that $f(v_*) \approx 2\cdot 10^{-17}$ for $d=5$ in~\eqref{boito811}); possibly for this reason, this was the only  numerical experiment on which Riemann-Oracle did not either match or outperform the (specialized) state of the art. It is possible that Riemann-Oracle can solve such badly scaled problems with minor improvements, for instance a more careful implementation of the optimization routine; we postpone this task to future research.

Another interesting comparison can be made with~\cite{UseM17gcd}. The matrix formulation of the GCD problem obtained in Section \ref{sec:approxgcd} essentially coincides with the one underlying the algorithm called $VP_S$ in~\cite{UseM17gcd}; however, \cite{UseM17gcd} then solves it by a different optimization method. The results of the $VP_S$ method reported in~\cite[Table~1]{UseM17gcd} match those in \cite{NAGASAKA} for $d=10,9,8$ only, while Riemann-Oracle performs better and also correctly solves $d=7,6$; see Table \ref{tbl:gcd}. The authors of \cite{UseM17gcd} also propose two additional formulations, labeled $VP_g$ and $VP_h$, which according to~\cite[Table~1]{UseM17gcd} perform better, although still not as well as \cite{NAGASAKA}. However, these formulations do not fall directly within our framework and thus are not (obviously) amenable to regularization as in Section~\ref{sec:need_regularization}. Potentially, this generates robustness issues; note that this particular experiment, which generates a full-rank problem in the sense of Theorem~\ref{thm:fv}, is not a reliable test in this sense. A possible question for future research is whether the regularization techniques described in the present paper can also be applied to the formulations of $VP_g$ and $VP_h$ in \cite{UseM17gcd}.

\subsection{Nearest matrix with prescribed nullity: A comparison with \cite{MarU14software}}\label{sub:prescribednullity}

We consider a structured matrix $A$, and aim to numerically approximate its distance to the nearest matrix with nullity at least $\ell$. We follow the strategy presented in Section \ref{sec:nearest_matrix_pescribed_nullity} and optimize over the Grassmann manifold \cite{EdelmanGrassmann}, by means of the augmented Lagrangian formulation. In this experiment, we consider the \texttt{grcar} matrix of size $8 \times 8$ 
\[
A = \begin{bmatrix}
     1 &  1 & 1  & 1  & 0  & 0  & 0 & 0\\
    -1 &  1 & 1 &  1 &  1 & 0 & 0 & 0\\
     0 & -1 & 1 & 1 & 1 & 1 & 0 & 0\\
     0 & 0 & -1 & 1 & 1 & 1 & 1 & 0\\
     0 & 0 & 0 & -1 & 1 & 1 & 1 & 1\\
     0 & 0 & 0 & 0 & -1 & 1 & 1 & 1\\
     0 & 0 & 0 & 0 & 0 & -1 & 1 & 1\\
     0 & 0 & 0 & 0 & 0 & 0 & -1 & 1
\end{bmatrix},
\]
with $\norm{A}_F=5.7446$. We choose a prescribed nullity $\ell \in \left\lbrace 1,2,3,4,5,6,7 \right\rbrace$, and compute
\begin{enumerate}[(a)]
\item the nearest matrix $\Delta$ with the same sparsity pattern of $A$, with nullity at least $\ell$;
\item the nearest Toeplitz matrix $\Delta$, with nullity at least $\ell$.
\end{enumerate}
We compare the output of our algorithm against the output of the software \cite{MarU14software}. For both algorithms, we used $[e_1 \ e_2 \ \cdots \ e_{\ell}]$ as a starting point for each ${\ell}$, i.e., the range of $[e_1 \ e_2 \ \cdots \ e_{\ell}]$ is the initial guess for the right kernel of the perturbed matrix.

In Table \ref{tab:grcar_nullity}, we report the results of the two experiments. In both numerical tests, Riemann-Oracle significantly outperforms the algorithm implemented in \cite{MarU14software}. Moreover, for experiment (b), the method in \cite{MarU14software} converges to the trivial solution $\Delta = -A$ for nullities 3 and higher. It is also worth noting that, for experiment (a), the method of this paper produces a non-decreasing sequence in the parameter ${\ell}$, while \cite{MarU14software} does not produce such a sequence. This suggests that the algorithm in \cite{MarU14software} may be more sensitive to initialization than our method.

\begin{table}[h!]
    \centering 
    \begin{tabular}{c  c c}
    \multicolumn{3}{c}{Sparsity pattern}\\
    \toprule
        $\ell$   & RO & \texttt{slra}\\
       \midrule 
        $1$ & $1.4126$ & 1.4142 \\
        $2$  & $2.1547$ & 2.6458 \\
        $3$  & $2.5905$ & 4.9615 \\
        $4$  & $3.2308$ & 3.6772 \\
        $5$  & $3.7762$ & 4.3589 \\
        $6$  & $4.4584$ & 5.7446 \\
        $7$ & $5.1418$ & 5.2915\\
       \bottomrule 
     \end{tabular}
     \hspace{1em}
      \begin{tabular}{c  c  c}
    \multicolumn{3}{c}{Toeplitz structure}\\
   \toprule 
        $\ell$   & RO & \texttt{slra}\\
        \midrule 
       $1$ & $1.2655$ & 1.5030\\
        $2$ & $1.8710$ & 5.7444\\
         $3$ & $2.2376$ & 5.7446\\
       $4$ & $3.0005$ & 5.7446\\
         $5$ & $3.3692$ & 5.7446\\
         $6$ & $4.1665$ & 5.7446\\
        $7$ & $5.0975$ & 5.7446 \\
       \bottomrule
     \end{tabular}
    \caption{ Left: Distance from $A$ to the nearest matrix with prescribed nullity and having the same zero pattern as $A$. Right: Distance from $A$ to the nearest Toeplitz matrix with prescribed nullity. RO refers to the Riemann-Oracle method of this paper, while \texttt{slra} refers to the variable projection method of \cite{MarU14software}.}
    \label{tab:grcar_nullity}
\end{table}

\section{Conclusions}\label{sec:conclusions}

We have introduced and studied Riemann-Oracle, an algorithm to solve matrix nearness problems that is applicable to a broad range of nearness problems of practical interest. The philosophy of Riemann-Oracle is based on the insight that some problems become much easier if an oracle reveals certain features of the minimizer, and one can then optimize over such features, typically over a real Riemannian manifold. We discussed the underlying theory and its links with least squares problems, and argued that it can lead to a powerful algorithm when enhanced with regularization tools from optimization. Extensive numerical experiments show that this method is already {competitive} on all test problems, and often provides a noticeable improvement with respect to previous methods. Possible future research directions include studying the performance of Riemann-Oracle on even more classes of nearness problems, as well as further improving its implementation details.

\section*{Acknowledgements}
We sincerely thank the authors of \cite{bora} for providing the Matlab code for their algorithm, and the authors of \cite{Sicilia} for providing some numerical results needed for the comparison in Section \ref{sec:betterthansicilia}. We also thank Froil\'{a}n Dopico, Etna Lindy, Konstantin Usevich, and two anonymous reviewers for providing valuable comments during the preparation of the manuscript. {MG and FP wish to thank Aalto University for hosting them during research visits that initiated this work. Furthermore, for part of the preparation of this work, MG was affiliated with Gran Sasso Science Institute, L'Aquila (Italy).}

\bibliographystyle{abbrv}
\bibliography{oracle}

\begin{thebibliography}{10}

\bibitem{AbsilBaker}
P.-A. Absil, C.~G. Baker, and K.~A. Gallivan.
\newblock Trust-region methods on {R}iemannian manifolds.
\newblock {\em Found. Comput. Math.}, 7(3):303–330, 2007.

\bibitem{AbsilMahony}
P.-A. Absil, R.~Mahony, and J.~Trumpf.
\newblock An extrinsic look at the {R}iemannian {H}essian.
\newblock In F.~Nielsen and F.~Barbaresco, editors, {\em Geometric Science of
  Information}, pages 361--368, Berlin, Heidelberg, 2013. Springer.

\bibitem{Bertsekas}
D.~Bertsekas.
\newblock {\em Nonlinear Programming}.
\newblock Athena Scientific, Belmont, MA, 1999.

\bibitem{Bini2010}
D.~A. Bini and P.~Boito.
\newblock A fast algorithm for approximate polynomial {GCD} based on structured
  matrix computations.
\newblock In D.~A. Bini, V.~Mehrmann, V.~Olshevsky, E.~E. Tyrtyshnikov, and
  M.~van Barel, editors, {\em Numerical Methods for Structured Matrices and
  Applications: The Georg Heinig Memorial Volume}. 2010.

\bibitem{Birgin}
E.~G. Birgin and J.~M. Martínez.
\newblock {\em Practical Augmented Lagrangian Methods for Constrained
  Optimization}.
\newblock SIAM, Philadelphia, PA, 2014.

\bibitem{Bjorck}
{\AA}.~Bj{\"o}rck.
\newblock {\em Numerical methods for least squares problems}.
\newblock SIAM, Philadelphia, PA, 1996.

\bibitem{borsdorf}
R.~Borsdorf.
\newblock An algorithm for finding an optimal projection of a symmetric matrix
  onto a diagonal matrix.
\newblock {\em SIAM J. Matrix Anal. Appl.}, 35(1):198--224, 2014.

\bibitem{Boumal}
N.~Boumal.
\newblock {\em An introduction to optimization on smooth manifolds}.
\newblock Cambridge University Press, Cambridge, 2023.

\bibitem{BoumalMishraAbsil}
N.~Boumal, B.~Mishra, P.-A. Absil, and R.~Sepulchre.
\newblock Manopt, a matlab toolbox for optimization on manifolds.
\newblock {\em Journal of Machine Learning Research}, 15(42):1455--1459, 2014.

\bibitem{Bye88}
R.~Byers.
\newblock A bisection method for measuring the distance of a stable matrix to
  the unstable matrices.
\newblock {\em SIAM J. Sci. Stat. Comput.}, 9(5):875--881, 1988.

\bibitem{nearsingpen}
R.~Byers, C.~He, and V.~Mehrmann.
\newblock Where is the nearest non-regular pencil?
\newblock {\em Linear Algebra Appl.}, 285(1):81--105, 1998.

\bibitem{bora}
B.~Das and S.~Bora.
\newblock Nearest rank deficient matrix polynomials.
\newblock {\em Linear Algebra Appl.}, 674:304--350, 2023.

\bibitem{dhillontropp}
I.~S. Dhillon and J.~A. Tropp.
\newblock Matrix nearness problems with {B}regman divergences.
\newblock {\em SIAM J. Matrix Anal. Appl.}, 29(4):1120--1146, 2007.

\bibitem{dd}
A.~Dmytryshyn and F.~M. Dopico.
\newblock Generic complete eigestructures for sets of matrix polynomials with
  bounded rank and degree.
\newblock {\em Linear Algebra Appl.}, 535:213--230, 2017.

\bibitem{DNN24}
F.~Dopico, V.~Noferini, and L.~Nyman.
\newblock A {R}iemannian optimization method to compute the nearest singular
  pencil.
\newblock {\em SIAM J. Matrix Anal. Appl.}, 45(4):2007--2038, 2024.

\bibitem{EY36}
C.~Eckart and G.~Young.
\newblock The approximation of one matrix by another of lower rank.
\newblock {\em Psychometrika}, 1(3):211--218, 1936.

\bibitem{EdelmanGrassmann}
A.~Edelman, T.~Arias, and S.~Smith.
\newblock The geometry of algorithms with orthogonality constraints.
\newblock {\em SIAM J. Matrix Anal. Appl.}, 20(2):303--353, 1998.

\bibitem{forney}
G.~D. Forney~Jr.
\newblock Minimal bases of rational vector spaces, with applications to
  multivariable linear systems.
\newblock {\em SIAM J. Control}, 13(3):493--520, 1975.

\bibitem{variableprojection}
G.~H. Golub and V.~Pereyra.
\newblock The differentiation of pseudo-inverses and nonlinear least squares
  problems whose variables separate.
\newblock {\em SIAM J. Numer. Anal.}, 10(2):413--432, 1973.

\bibitem{GCE23}
F.~Goyens, C.~Cartis, and A.~Eftekhari.
\newblock Nonlinear matrix recovery using optimization on the {G}rassmann
  manifold.
\newblock {\em Appl. Comput. Harmon. Anal.}, 62:498–542, 2023.

\bibitem{Grabner}
N.~Gräbner, V.~Mehrmann, S.~Quraishi, C.~Schröder, and U.~von Wagner.
\newblock Numerical methods for parametric model reduction in the simulation of
  disk brake squeal.
\newblock {\em ZAMM - Journal of Applied Mathematics and Mechanics /
  Zeitschrift für Angewandte Mathematik und Mechanik}, 96(12):1388--1405,
  2016.

\bibitem{Sicilia}
N.~Guglielmi, C.~Lubich, and S.~Sicilia.
\newblock Rank-1 matrix differential equations for structured eigenvalue
  optimization.
\newblock {\em SIAM J. Numer. Anal.}, 61(4):1737--1762, 2023.

\bibitem{HeWatson}
C.~He and G.~A. Watson.
\newblock An algorithm for computing the distance to instability.
\newblock {\em SIAM J. Matrix Anal. Appl.}, 20(1):101--116, 1998.

\bibitem{nickthesis}
N.~J. Higham.
\newblock {\em Nearness Problems in Numerical Linear Algebra}.
\newblock PhD thesis, University of Manchester, 1985.

\bibitem{nickspd}
N.~J. Higham.
\newblock Computing a nearest symmetric positive semidefinite matrix.
\newblock {\em Linear Algebra Appl.}, 103:103--118, 1988.

\bibitem{nicksurvey}
N.~J. Higham.
\newblock Matrix nearness problems and applications.
\newblock In {\em M. J. C. Gover and S. Barnett, editors, Applications of
  Matrix Theory}, pages 1--27. Oxford University Press, Oxford, UK, 1989.

\bibitem{nickcorr}
N.~J. Higham.
\newblock Computing the nearest correlation matrix—a problem from finance.
\newblock {\em IMA J. Numer. Anal.}, 22(3):329--343, 2002.

\bibitem{KYZ}
E.~Kaltofen, Z.~Yang, and L.~Zhi.
\newblock Structured low rank approximation of a {S}ylvester matrix.
\newblock In D.~Wang and L.~Zhi, editors, {\em Symbolic-Numeric Computation},
  pages 69--83, Basel, 2007. Birkh{\"a}user Basel.

\bibitem{KOM09}
R.~H. Keshavan, A.~Montanari, and S.~Oh.
\newblock Matrix completion from a few entries.
\newblock {\em IEEE Transactions on Information Theory}, 56(6):2980--2998,
  2010.

\bibitem{Kre06}
D.~Kressner.
\newblock Finding the distance to instability of a large sparse matrix.
\newblock In {\em 2006 IEEE Conference on Computer Aided Control System Design,
  2006 IEEE International Conference on Control Applications, 2006 IEEE
  International Symposium on Intelligent Control}, pages 31--35, 2006.

\bibitem{LiuBoumal}
C.~Liu and N.~Boumal.
\newblock Simple algorithms for optimization on {R}iemannian manifolds with
  constraints.
\newblock {\em Appl. Math. Optim.}, 82(3):949–981, 2020.

\bibitem{MarU13missing}
I.~Markovsky and K.~Usevich.
\newblock Structured low-rank approximation with missing data.
\newblock {\em SIAM J. Matrix Anal. Appl.}, 34(2):814--830, 2013.

\bibitem{MarU14software}
I.~Markovsky and K.~Usevich.
\newblock Software for weighted structured low-rank approximation.
\newblock {\em J. Comput. Appl. Math.}, 256:278--292, 2014.

\bibitem{NAGASAKA}
K.~Nagasaka.
\newblock Toward the best algorithm for approximate {G}{C}{D} of univariate
  polynomials.
\newblock {\em J. Symbolic Comput.}, 105:4--27, 2021.

\bibitem{NoceWrig}
J.~Nocedal and S.~J. Wright.
\newblock {\em Numerical Optimization}.
\newblock Springer, New York, NY, USA, 2nd edition, 2006.

\bibitem{noferini24}
V.~Noferini.
\newblock Invertible bases and root vectors for analytic matrix-valued
  functions.
\newblock {\em Electron. J. Linear Algebra}, 40:1--13, 2024.

\bibitem{NP21}
V.~Noferini and F.~Poloni.
\newblock Nearest {$\Omega$}-stable matrix via {R}iemannian optimization.
\newblock {\em Numer. Math.}, 148(4):817--851, 2021.

\bibitem{oleary}
D.~O’Leary and B.~Rust.
\newblock Variable projection for nonlinear least squares problems.
\newblock {\em Comput. Optim. Appl.}, 54(3):579–593, 2013.

\bibitem{papa}
C.~H. Papadimitriou.
\newblock {\em Computational {C}omplexity.}
\newblock Addison-Wesley, 1994.

\bibitem{qisun}
H.~Qi and D.~Sun.
\newblock A quadratically convergent {N}ewton method for computing the nearest
  correlation matrix.
\newblock {\em SIAM J. Matrix Anal. Appl.}, 28(2):360--385, 2006.

\bibitem{RPG}
J.~B. Rosen, H.~Park, and J.~Glick.
\newblock Total least norm formulation and solution for structured problems.
\newblock {\em SIAM J. Matrix Anal. Appl.}, 17(1):110--126, 1996.

\bibitem{sipser2006}
M.~Sipser.
\newblock {\em Introduction to the Theory of Computation}.
\newblock Course Technology, Boston, MA, 2nd edition, 2006.

\bibitem{TM}
F.~Tisseur and K.~Meerbergen.
\newblock The quadratic eigenvalue problem.
\newblock {\em SIAM Review}, 43(2):235--286, 2001.

\bibitem{UseM14manifold}
K.~Usevich and I.~Markovsky.
\newblock Optimization on a {Grassmann} manifold with application to system
  identification.
\newblock {\em Automatica}, 50(6):1656--1662, 2014.

\bibitem{UseM14}
K.~Usevich and I.~Markovsky.
\newblock Variable projection for affinely structured low-rank approximation in
  weighted {{\(2\)}}-norms.
\newblock {\em J. Comput. Appl. Math.}, 272:430--448, 2014.

\bibitem{UseM17gcd}
K.~Usevich and I.~Markovsky.
\newblock Variable projection methods for approximate (greatest) common divisor
  computations.
\newblock {\em Theor. Comput. Sci.}, 681:176--198, 2017.

\bibitem{VL85}
C.~Van~Loan.
\newblock How near is a stable matrix to an unstable matrix?
\newblock In {\em Linear algebra and its role in systems theory ({B}runswick,
  {M}aine, 1984)}, volume~47 of {\em Contemp. Math.}, pages 465--478. Amer.
  Math. Soc., Providence, RI, 1985.

\bibitem{bartvander}
B.~Vandereycken.
\newblock Low-rank matrix completion by {R}iemannian optimization.
\newblock {\em SIAM Journal on Optimization}, 23(2):1214--1236, 2013.

\bibitem{ist}
G.~Verghese, P.~Van~Dooren, and T.~Kailath.
\newblock Properties of the system matrix of a generalized state-space system.
\newblock {\em Internat. J. Control}, 30(2):235--243, 1979.

\bibitem{Zhang}
W.~H. Yang, L.-H. Zhang, and R.~Song.
\newblock Optimality conditions for the nonlinear programming problems on
  {R}iemannian manifolds.
\newblock {\em Pacific Journal of Optimization}, 10(2):415--434, 2014.

\bibitem{Zeng}
Z.~Zeng.
\newblock The numerical greatest common divisor of univariate polynomials.
\newblock In {\em Randomization, relaxation, and complexity in polynomial
  equation solving}, volume 556 of {\em Contemp. Math.}, pages 187--217.
  American Mathematical Society, Providence, RI, 2011.

\end{thebibliography}

\appendix
\section{Proof of Theorem~\ref{thm:gradienthessian}} \label{sec:proofgradienthessian}
To determine the gradient $\nabla f_{\varepsilon,y}(v)$ of $f_{\varepsilon,y}(v)$, we compute the directional derivative in a direction $w$
\[
\dot{f}_{\varepsilon,y} = \frac{d}{d t} f_{\varepsilon,y}(v+tw)\bigg|_{t=0}.
\]
We use a dotted letter to denote the derivative (in $t$) of the corresponding non-dotted quantity, evaluated in $t=0$; hence for instance
\[
\dot{r} = \frac{d}{d t} r(v+tw)\bigg|_{t=0} = -Aw.
\]
It is clear that the formulae in~\eqref{dotstuff} hold, when ascribing this meaning to dotted letters. We now show that~\eqref{dotz} holds. Using the expansion $(X+H)^{-1} = X^{-1} - X^{-1}HX^{-1} + O(\norm{H}^2)$ with $X = MM^* + \varepsilon I$, we have
\[
\dot{z} = X^{-1}r = -X^{-1}\dot{X}X^{-1}r + X^{-1}\dot{r} = -X^{-1}(\dot{M}M^* + M\dot{M}^*)z - X^{-1}Aw.
\]
We now observe that 
\[
    \dot{M}M^*z = \dot{M}\delta_* = \sum_{i=1}^p P^{(i)}w (\delta_*)_i = \Delta_* w
\]
to obtain the formula~\eqref{dotz}.

We can now compute the derivative
\begin{align*}
\dot{f}_{\varepsilon,y} &= \dot{r}^*z + r^*\dot{z}\\
&= -w^*A^* z - \underbrace{r^* (MM^*+\varepsilon I)^{-1}}_{=z^*}(M\dot{M}^*z + (A+\Delta_*)w)\\
&= -w^*A^*z - z^*Aw - \underbrace{z^* M\dot{M}^*}_{=(\Delta_*w)^*}z  - z^*\Delta_*w\\
&= -2\operatorname{Re} z^*(A+\Delta_*)w\\
&= \left\langle -2(A+\Delta_*)^*z, w \right\rangle.
\end{align*}
This equality, which holds for all $w$, proves the formula for the gradient~\eqref{gradf}.

To compute the action of the Hessian, we differentiate the gradient to obtain
\begin{align*}
    H_{\varepsilon,y}(v) w &= \frac{d}{dt} \nabla f_{\varepsilon,y}(v+tw)\bigg|_{t=0} \\&= -2(\underbrace{\dot{A}}_{=0}+\dot{\Delta}_*)^*z -2 (A+\Delta_*)^*\dot{z}. \qedhere
\end{align*}

\begin{remark}
    It is interesting to note the alternative form
    \begin{equation*}
        (\Delta_*)^*z = \left(\sum_{i=1}^p P^{(i)}((P^{(i)}v)^*z))\right)^* z = \left(\sum_{i=1}^p P^{(i)}zz^*(P^{(i)})^* \right)v;
    \end{equation*}
    which shows that when $y=0$ the gradient can be written as the product of a Hermitian matrix with $v$:
    \[
        \nabla f_{\varepsilon,0}(v) = 2\left(A^*(MM^*+\varepsilon I)^{-1}A - \sum_{i=1}^p P^{(i)}zz^*(P^{(i)})^*\right)v.
    \]
\end{remark}

\end{document}